\documentclass[graybox]{svmult}

\usepackage{mathptmx}       
\usepackage{helvet}         
\usepackage{courier}        
\usepackage{type1cm}        
\usepackage{makeidx}         
\usepackage{graphicx}        
\usepackage{multicol}        
\usepackage[bottom]{footmisc}

\usepackage{newtxtext,newtxmath}
\usepackage[numbers]{natbib}
\usepackage{url}
\usepackage{xcolor}
\definecolor{newcolor}{rgb}{.8,.349,.1}

\usepackage{comment}
\usepackage{newtxtext,newtxmath}
\usepackage{latexsym}

\DeclareMathOperator{\atantwo}{atan2}
\usepackage{hyperref}

\usepackage{mathtools}
\mathtoolsset{showonlyrefs=true}
\usepackage[most]{tcolorbox}
\usepackage{graphicx}
\usepackage{tikz}
 \usepackage{geometry}
 \newgeometry{vmargin={25mm}, hmargin={25mm,25mm}} 
\usepackage{amsmath}

\makeindex             

\graphicspath{{WilliamsFigures/}}

\begin{document}

\title*{The architectural application of shells whose boundaries subtend a constant solid angle}
\author{Emil Adiels, Mats Ander  
and Chris J. K. Williams}
\institute{Emil Adiels \at Chalmers University of Technology, Sweden, \email{emil.adiels@chalmers.se}
\and Mats Ander \at Chalmers University of Technology, Sweden, \email{mats.ander@chalmers.se}
\and Chris J. K. Williams \at Chalmers University of Technology, Sweden, \email{christopher.williams@chalmers.se}}
%
%
\maketitle

\abstract{ 
Surface geometry plays a central role in the design of bridges, vaults and shells, using various techniques for generating a geometry which aims to balance structural, spatial, aesthetic and construction requirements.
\newline \indent
In this paper we propose the use of surfaces defined such that given closed curves subtend a constant solid angle at all points on the surface and form its boundary. Constant solid angle surfaces enable one to control the boundary slope and hence achieve an approximately constant span-to-height ratio as the span varies, making them structurally viable for shell structures. In addition, when the entire surface boundary is in the same plane, the slope of the surface around the boundary is constant and thus follows a principal curvature direction. Such surfaces are suitable for surface grids where planar quadrilaterals meet the surface boundaries. They can also be used as the Airy stress function in the form finding of shells having forces concentrated at the corners.
\newline \indent
Our technique employs the Gauss-Bonnet theorem to calculate the solid angle of a point in space and Newton's method to move the point onto the constant solid angle surface. We use the Biot-Savart law to find the gradient of the solid angle. The technique can be applied in parallel to each surface point without an initial mesh, opening up for future studies and other applications when boundary curves are known but the initial topology is unknown.
\newline \indent
We show the geometrical properties, possibilities and limitations of surfaces of constant solid angle using examples in three dimensions.}

\section{Introduction}\label{sec1}
\begin{figure}[t]
\centering
\includegraphics[width=1.0\textwidth]{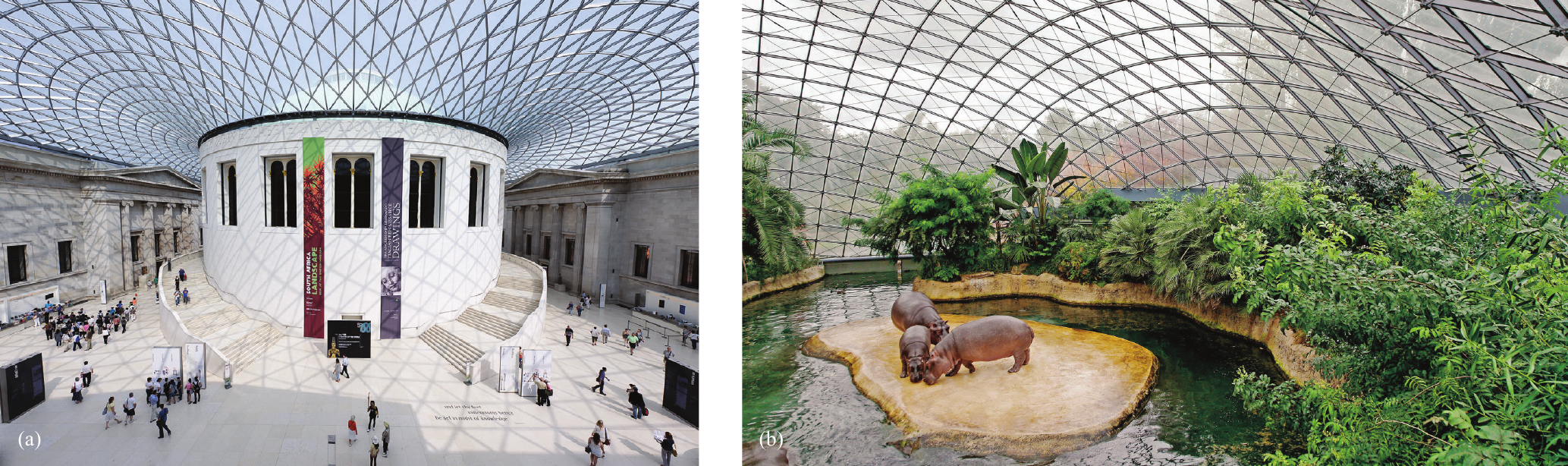}
\caption{ (a) The British Museum Great Court roof, by Foster and partners, Buro Happold and Waagner Biro. Photo by Andrew Stawarz, licensed under CC BY-ND 2.0. (b) grid shell of the Hippo House by Schlaich Bergermann partner, photo used with permission © Zoo Berlin }
\label{fig:britishCourt}
\end{figure}

Architectural geometry is the application of geometry to the design and construction of buildings and bridges, particularly those with curved surfaces like shells and grid shells (Figs. \ref{fig:britishCourt}(a) and (b)). The surface curvature enables the shell to carry load mainly through membrane action, making them much more material efficient than conventional flat slabs and beams used today. The complex geometry, combined with requirements and desires regarding economic, structural, production, spatial and aesthetic aspects, makes this a topic that has fascinated builders and mathematicians for centuries. Early treatises in architectural geometry include \textit{Le premier tome de l'architecture} by Philibert de Lorme (1512-1570), examining the art of cutting stones in vaults, while an extensive overview of contemporary techniques and applications from the field of differential geometry can be found in Pottmann et al. \cite{Pottmann2014}. Meanwhile, shell designers have experimented with various  shapes to balance requirements and qualities, like ruled surfaces by Antoni Gaudí \citep{Collins1963,Huerta2006} and Félix Candela \citep{Faber1963}, Eladio Dieste’s "Gaussian vaults" \citep{Allen}, and translational surfaces (Fig. \ref{fig:britishCourt}(b)) by Jörg Schlaich \citep{Schlaich1996}. Other examples include Weingarten surfaces \cite{Pellis2021}, such as minimal surfaces, surfaces of revolution and constant mean surfaces. Additional techniques include form finding \cite{Adriaenssens2014} striving for structural efficiency or a specific state of stress for a given load. Examples include minimising surface tension, bending energy or strain energy producing soap films \cite{Klaus1987}, Willmore energy surfaces \cite{WILLIAMS1987310} or hanging nets \cite{RubioiBellver1913}. 

\begin{figure}[b!]
\centering
\includegraphics[width=1.0\textwidth]{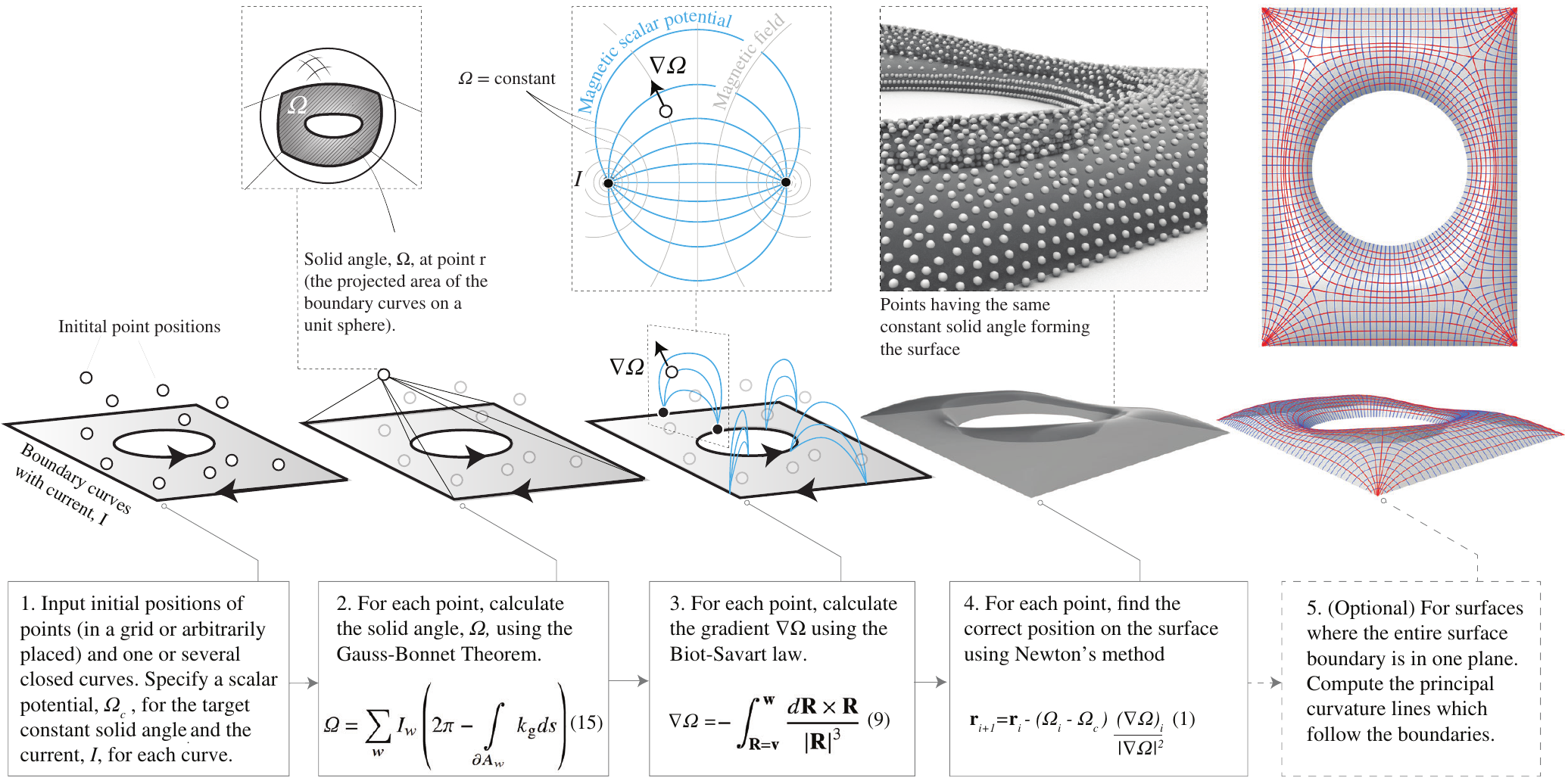}
\caption{The steps of our method for generating surfaces of constant solid angle. The first four steps are necessary while the last step is optional and specifically interesting for surfaces where the entire boundary is in the same plane.}
\label{fig:technique}
\end{figure}

We propose an alternative using surfaces defined by points, all having the same constant solid angle, $\mathit{\Omega_c}$, subtended by its boundary curves. The mathematical problem we aim to solve for generating these surfaces can be formulated as the solution to \eqref{eq:newton1} approximated using Newton's method (step 4 in Fig. \ref{fig:technique}) 

\begin{equation}\label{eq:newton1}
    \mathbf{r}_{i+1}=\mathbf{r}_i-\left(\mathit{\Omega}_i-\mathit{\Omega}_{\mathrm{c}}\right)\frac{\left(\nabla\mathit{\Omega}\right)_i}{\left|\nabla\mathit{\Omega}\right|_i^2} \tag{1}\mathrm{,} 
\end{equation}
 for each point \textbf{r}, which includes the task of formulating the expressions for computing the solid angle $\mathit{\Omega}_i$ and the gradient $(\nabla\mathit{\Omega})_i$ for an arbitrary point in space. 

The benefit of constant solid angle surfaces is that they enable one to control the boundary slope and hence achieve an approximately constant span-to-height ratio as the span varies. Thus, making them structurally viable for shells and grid shells, even though the formulation takes no regard of structural aspects. For cases where the entire boundary is in the same plane, the slope is constant along its boundaries. This means the surface boundary follows a principal curvature direction, which is an otherwise rare occurrence. Therefore, suitable for situations when surface grids of planar quadrilaterals are desired for economic or production requirements since the panels are not cut at the boundary, which usually is an unfortunate consequence of these patterns as seen in Fig. \ref{fig:britishCourt}(b). Another potential application for surfaces having constant slope along the surface boundary is as Airy stress functions, $ \phi$, which can be used in Pucher's equation \cite{Timoshenko1959} for form finding of shells as done by Miki. et al \cite{Miki2015, Miki2022}. Such an Airy stress function would result in a shell where the lateral thrust is concentrated to the corners. Thus, suitable for situations similar to that of British Museum Great Court roof where the roof lies on roller supports such that the lateral thrust needs to be handled in the corners.

Our technique for generating these surfaces is simple and straightforward since it does not necessarily require any initial mesh, such as many form finding techniques used in architecture, for example dynamic relaxation \cite{Adriaenssens2014} or force density method \cite{Schek1974}. Each point can find itself onto the surface independently meaning the process can be done in parallel on a GPU. It opens up a third application for this type of surfaces beyond architecture in for instance industrial design or computer graphics. Having only known boundaries one can find a surface using this technique without having an initial mesh. 

The method, illustrated in Fig. \ref{fig:technique}, at minimum, requires the designer to input one closed boundary curve with a chosen strength and direction of the current, a scalar potential for the constant solid angle, and a number of points that can be positioned in a grid or arbitrarily placed. Based on the inputs from the designer there are then three more steps. First, we calculate the solid angle using the Gauss-Bonnet theorem, section \ref{solidAngleCalculation}.
Secondly, we find the gradient, which is also the surface normal, using the Biot-Savart law, section \ref{solidAngleGradient}. 
Lastly, we use Newton's method to move the points to the surface, section \ref{sec:newton}. Examples of various surfaces can be found in section \ref{sec:examples}. The technique we present is for generating the surface itself, which is the main focus of this paper. However, it could be expanded to include the generation of surface grids such as principal curvature nets, a Chebyshev nets, and geodesic coordinates. Optionally, such grids can be applied afterwards on the generated surface. In section \ref{sec:curvature} we show two methods to compute the curvature and principal curvature directions, which can be used to generate a principal curvature net.

\subsection{Connection between the shell  form and the shell pattern}
The surface and its surface pattern of triangles or quadrilaterals is connected through the components of the first and the second fundamental form \citep{Struik1961, Stoker1969, GreenZerna68}. From a structural point of view, a triangular grid might be more desirable than a quadrilateral grid, but from a production point of view, a quadrilateral grid may be preferable as described by \citeauthor{Schlaich2005} \cite{Schlaich2005}.

There are three compatibility equations, the Gauss-Codazzi equations, connecting the six components of the first and second fundamental form, which tell us how a pattern on a surface has to be deformed to cover a doubly curved surface. If the surface can be considered to be acting as a membrane shell there are three components of the stress tensor and three equations of static equilibrium \citep{GreenZerna68}. Meaning there are, in total, nine unknowns and six equations requiring the designer to introduce three more equations or constraints, which could be used to facilitate a more economic production. Examples include equal mesh Chebyshev nets\footnote{ Chebyshev nets require constraints $a_{11} = a_{22} = L, n^{12} =0$ using notation in \citeauthor{GreenZerna68} \cite{GreenZerna68}} \citep{Chebyshev1946}  used for continuous laths in the Mannheim Multihalle \citep{Liddell2015}, geodesic coordinate nets\footnote{Geodesic coordinate nets require constraints $a_{22} = L, a_{12} = n^{12} = 0 $ using notation in  \citeauthor{GreenZerna68} \cite{GreenZerna68}} following the bed joints on masonry shells \citep{Adiels2017,Adiels2021} or cutting patterns of tensile nets \citep{Williams1980}, and principal curvature nets.

Principal curvature nets consist of principal curvature lines intersecting at right angles. From a production point of view, it means approximately torsion-free nodes and panels that can be made approximately planar, at least if the grid is fine. If the grid is not fine then one is in the realm of discrete differential geometry \citep{Crane:2013:DGP,10.1007/978-3-319-66272-5_1}. In terms of the first and second fundamental form of a continuous surface it means \textit{F} and \textit{f} are zero using the notation in Struik \cite{Struik1961}($a_{12}$ and $b_{12}$ in Green \& Zerna \cite{GreenZerna68}). Work has also been done in aligning the principal stress and principal curvature direction by for instance Tellier et al. \cite{Tellier2019} and Pellis et al. \cite{Pellis2018}. One of the issues with a principal curvature net is that it does not guarantee a nice connection to the boundary, making it necessary to cut the grid, which is usually not good for the architectural expression and requires specially made panels and components. Thus, ideally, one would have a surface such that the principal directions are aligned with the boundaries, meaning a constant slope along the boundary.
\subsection{Surfaces exercising boundary tangency restrictions}
Previous work in controlling the tangency constraints along the boundary have been done using surfaces which minimize the Willmore energy \citep{Bobenko2005}. This is equivalent to minimizing the bending energy and is similar to the shapes we see in cells \citep{Muller2014}, but it is also a feature that is beneficial for structural shells and membranes \citep{WILLIAMS1987310}. 

To find a surface which minimizes the Willmore energy it is necessary to solve a differential equation whereas points on a surface of constant solid angle subtended by the boundary curves but can be found individually and independently for each point. When the surface boundaries are in the same plane the slope along the boundary is constant, a constant one can choose. Surface of constant solid angle appear in potential theory \citep{Lamb1932} and can be used for very complicated boundaries \citep{Binysh2019,Binysh2018}.

A simpler version of the technique described in this paper was used for the design of the courtyard roofs of the British Museum Great Court \citep{Williams2001} (Fig. \ref{fig:britishCourt}(a)) and Het Scheepvaartmuseum in Amsterdam \citep{Adriaenssens2012}. This simpler version used the reciprocal of the gradient of the solid angle in the horizontal plane to determine the vertical coordinate. In the case of a boundary consisting of two infinite parallel straight lines this would produce a parabolic cross section, rather than the curves shown in Fig. \ref{fig:gradientSolidAngle}.
\section{Surfaces of constant solid angle}
\label{sec:surfaceDef}
\begin{figure}[h]
\centering
\includegraphics[width=0.7\textwidth]{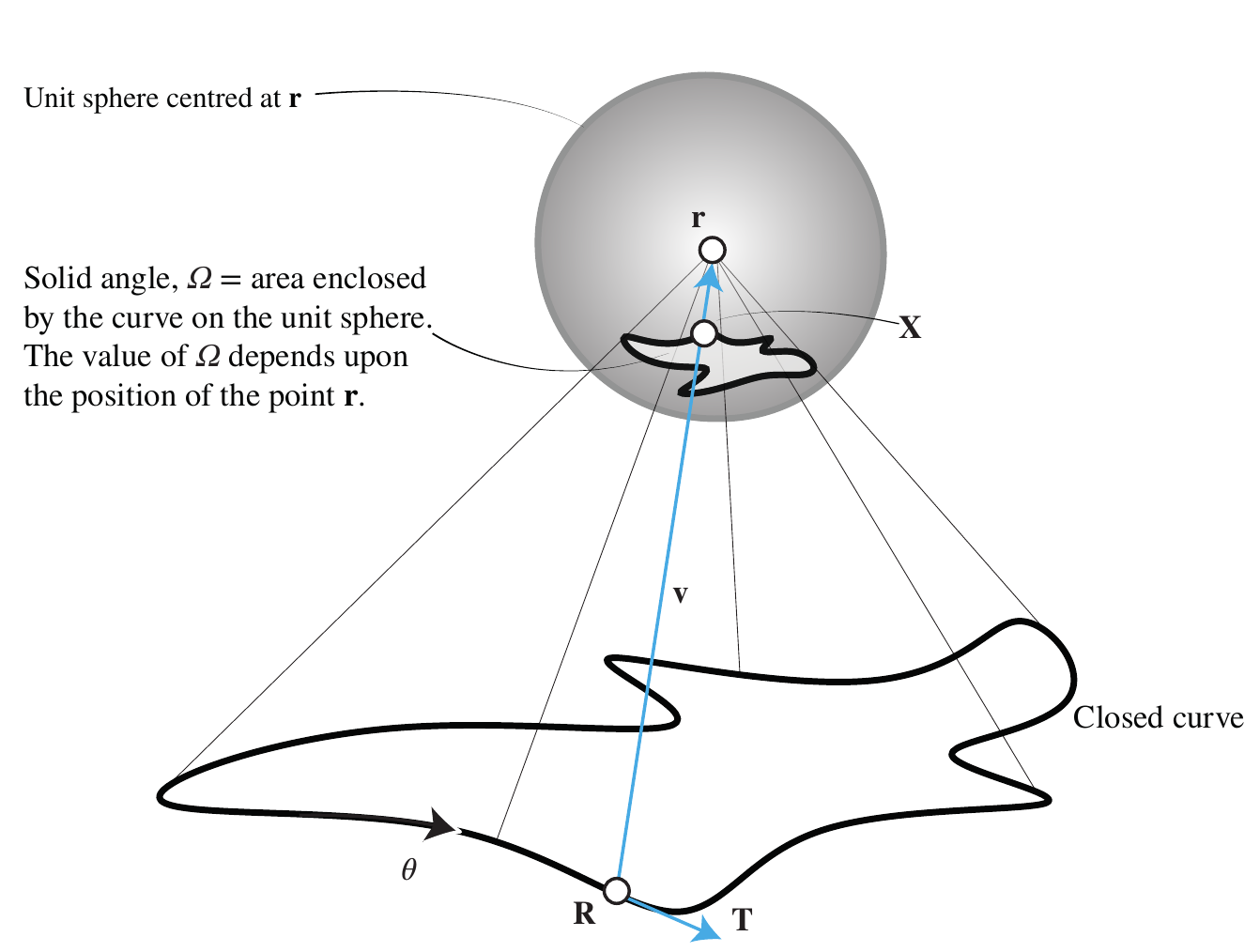}
\caption{Solid angle  $\mathit{\Omega}$ subtended at a point $\mathbf{r}$ is the surface area enclosed
by the curve on the unit sphere. $\mathbf{R}$ is a typical point on the boundary, which is a function of curve parameter $\theta$, where $\mathbf{T}$ is the tangent. Vector $\mathbf{v}$ goes from $\mathbf{R}$ to $\mathbf{r}$ and cuts the unit sphere at the point $\mathbf{X}$.
\newline \indent
$\mathbf{R}$, $\mathbf{r}$ and $\mathbf{X}$ are position vectors relative to some fixed origin. }
\label{fig:SolidAngle}
\end{figure}

An arbitrary point in space can be joined with straight lines to all the points on a given closed curve to form a ruled surface or cone, but almost invariably without a circular cross-section. This cone will intersect a sphere of arbitrary radius with another closed curve and the solid angle, measured in steradians, is the surface area enclosed by the curve on the sphere divided by the square of the radius of the sphere. The sphere is often taken as unit radius, but even so it is important to note that the solid angle is always dimensionless (Fig. \ref{fig:SolidAngle}).

\begin{figure}[htp]
\centering
\includegraphics[width=0.7\textwidth]{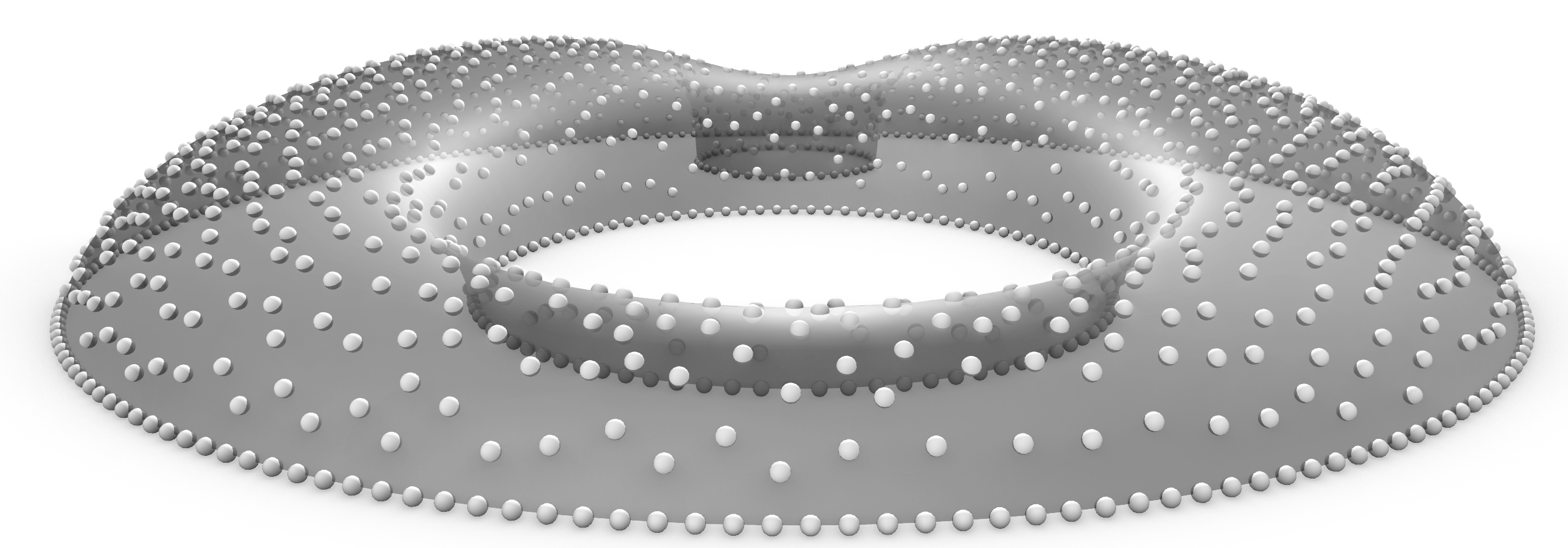}
\caption{Points having the same constant solid angle subtended by three circles in the same plane. They all lie on the constant solid angle surface.}
\label{fig:solidangleSurface}
\end{figure}

If the curve crosses itself on the unit sphere some areas will be positive and some negative, exactly as in ordinary integration.
The surface area of a complete unit sphere is $4\pi$ and if $\mathit{\Omega}$ is the solid angle subtended by the boundary curve, then $4\pi - \mathit{\Omega}$ is the area on the surface outside the boundary. In the case of a plane boundary, the solid angle is $2\pi$ if the apex is in the plane inside the boundary or 0 (or $4\pi$) if it is in the plane outside the boundary.
Going around any closed boundary a number of times the solid angle increases or decreases by $4\pi$ each revolution and then jumps back again. Being very close to a boundary the change in solid angle is twice the angle of rotation around the boundary, $4\pi$ instead of $2\pi$ for a complete rotation.
\newline \indent
It is now possible to define a \emph{constant solid angle surface} as the locus of points such that one or several given closed curves subtend the same solid angle $\mathit{\Omega}$ at all points on the surface, as in Fig. \ref{fig:solidangleSurface}. In other words the solid angle subtended by a given boundary curve at point $\mathbf{r}$ depends upon the position of $\mathbf{r}$. If we say that the solid angle should remain constant then $\mathbf{r}$ is constrained to move on a surface. Thus the surface is defined by a potentially infinite number of points on the surface and there are no $\textit{u,v}$ surface coordinates unless we constrain the points to lie on some grid. To define a particular point $\mathbf{r}$ we need to specify where it is on the surface, for example by specifying its position in plan.
\begin{tcolorbox}[colback=white,colframe=red!75!black]
  For a plane horizontal boundary, the slope of a constant solid angle surface at the boundary is constant and equal to $\pi - \dfrac{\mathit{\mathit{\Omega}}}{2}$. Thus, as we move around the boundary there is no rotation of the normal about the tangent to the boundary. The boundary is therefore a principal curvature direction on the surface. This is a special case of Joachimsthal's theorem \citep{Joachimsthal1846}.
\end{tcolorbox}

\section{Method of finding a constant solid angle surface} \label{sec:newton}
It is not possible to explicitly calculate the shape of a constant solid angle surface. However, we shall see that it is possible to calculate the solid angle $\mathit{\Omega}$ subtended by the boundary at any given point in space. We can also calculate the gradient of the solid angle $\nabla\mathit{\Omega}$ which tells us how the solid angle varies if we move slightly.

Let us imagine that $\mathit{\Omega}_{\mathrm{c}}$ is the constant value of solid angle on the surface we wish to generate. If we start at an arbitrary point in space $\mathbf{r}_i$ we can calculate the solid angle $\mathit{\Omega}_i$ and its gradient $\left(\nabla\mathit{\Omega}\right)_i$ at that point.

\begin{figure}[htp]
\centering
\includegraphics[width=0.7\textwidth]{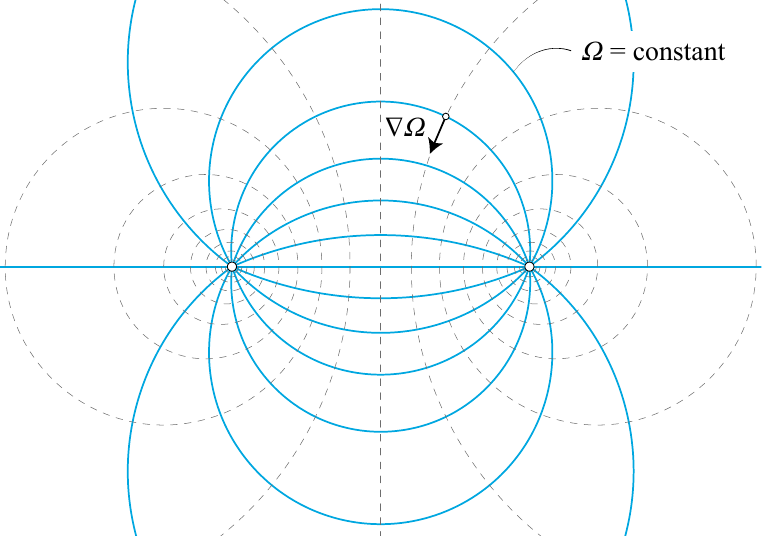}
\caption{Section through different level surfaces of constant solid angle and their gradient for the case of a boundary consisting of two infinitely long parallel straight lines. The surfaces correspond to surfaces of constant magnetic scalar potential for the case of a current running along one of the wires and back along the other. They also correspond to the surfaces of constant velocity potential associated with a clockwise and anticlockwise vortex similar to that left behind an aircraft.}
\label{fig:gradientSolidAngle}
\end{figure}

We can now move to a new point $\mathbf{r}_{i+1}$
\begin{equation}\label{eq:newtonCHRIS}
    \mathbf{r}_{i+1}=\mathbf{r}_i-\left(\mathit{\Omega}_i-\mathit{\Omega}_{\mathrm{c}}\right)\frac{\left(\nabla\mathit{\Omega}\right)_i}{\left|\nabla\mathit{\Omega}\right|_i^2}
\end{equation}
which will be nearer to the surface. We can expect to have to do this a number of times, but once we are near the surface it will converge rapidly as shown in Fig. \ref{fig:gradientSolidAngle}. This is an application of Newton's method, which is usually written
\begin{equation}
    x_{i+1}=x_i-\frac{f\left(x_i\right)}{f'\left(x_i\right)}
\end{equation}
for the solution of
\begin{equation}
    f\left(x\right)=0\mathrm{.}
\end{equation}
As always with Newton's method we have to be careful when the gradient is small in case we jump much too far, in which case we can multiply the movement in \eqref{eq:newtonCHRIS} by some factor less than 1.0.

\begin{tcolorbox}[colback=white,colframe=red!75!black]
The movement in \eqref{eq:newtonCHRIS} is in the direction of the gradient, normal to  a surface of constant $\mathit{\Omega}$. We are, of course, at liberty to also move points tangential to the surface in order to satisfy some requirement of the grid pattern on the surface.
\end{tcolorbox}
\section{Calculation of the solid angle subtended by a closed curve}
\label{solidAngleCalculation}
In order to find points on this solid angle surface it is necessary be able to calculate the solid angle subtended by the boundary at an arbitrary point in space and the gradient of the solid angle so that we can move points bit by bit onto the solid angle surface.
In order to find the solid angle it is necessary to do a double integral to find the surface area enclosed by the curve on the unit sphere. However, using the Gauss-Bonnet theorem \citep{Struik1961,Eisenhart1947} it is possible to reduce the double integral to a line integral around the curve,
\begin{equation*}
\int\limits_A {KdA}  + \oint\limits_{\partial A} {{k_{\mathrm{g}}}ds}  = 2\pi
\end{equation*}
in which $A$ is the area on the surface, which here is the unit sphere. $\partial A$ is its boundary with arc length $s$, $K$ is the Gaussian curvature of the surface, which is 1 on the unit sphere and $k_{\mathrm{g}}$ is the geodesic curvature of the boundary curve on the unit sphere. The geodesic curvature of a curve on a surface is the curvature of the curve as seen when looking at it directly back down the normal to the surface.

Note that the sign of $k_{\mathrm{g}}$ is changed if the direction of travel is reversed around the curve, with a corresponding change to the value of $\mathit{\Omega}$. The direction of travel does not matter as long as it is consistent. However, having a boundary consisting of several closed boundary curves the directions of travel must be coordinated and compatible.

Thus
\begin{equation}
\label{Omega}
\mathit{\Omega}  = 2\pi  - \oint\limits_{\partial A} {{k_{\mathrm{g}}}ds}\mathrm{.}
\end{equation}
In Fig. \ref{fig:SolidAngle} $\mathbf{R}\left(\theta\right)$ is a typical point on the boundary, which is a function of curve parameter $\theta$ which may or may not be equal to the arc length. $\mathbf{r}$ is the point in space at which we want to find the solid angle subtended by the boundary and $\mathbf{v}$ is the vector from $\mathbf{R}$ to $\mathbf{r}$.

Vector $\mathbf{v}$ cuts the unit sphere at the point given by the vector $\mathbf{X}$ and $\mathbf{T}$ is the unit tangent to the boundary. Thus we have
\begin{equation}
    \begin{split}
        \mathbf{R}_\theta&=\frac{d\mathbf{R}}{d\theta}\\
        \mathbf{T}&=\frac{\mathbf{R}_\theta}{\left|\mathbf{R}_\theta\right|}\\
        \mathbf{X}_\theta&=\frac{\left|\mathbf{R}_\theta\right|}{\left|\mathbf{v}\right|}\left(\mathbf{T}
        -\frac{\left(\mathbf{v}\cdot\mathbf{T}\right)\mathbf{v}}{\left|\mathbf{v}\right|^2}\right)\\
        \left|\mathbf{X}_\theta\right|&=\frac{\left|\mathbf{R}_\theta\right|\left|\mathbf{v}\times\mathbf{T}\right|}{\left|\mathbf{v}\right|^2}
    \end{split}
\end{equation}
and the geodesic curvature of the curve on the sphere is equal to
\begin{equation}
    \begin{split}
        k_{\mathrm{g}}&=\frac{1}{\left|\mathbf{X}_\theta\right|}\frac{d}{d\theta}\left(\frac{\mathbf{X}_\theta}{\left|\mathbf{X}_\theta\right|}\right)
        \cdot\frac{\left(\mathbf{v}\times\mathbf{T}\right)}{\left|\mathbf{v}\times\mathbf{T}\right|}
        =\frac{\left|\mathbf{v}\right|\left(\mathbf{T}_\theta\cdot\left(\mathbf{v}\times\mathbf{T}\right)\right)}{\left|\mathbf{X}_\theta\right|\left|\mathbf{v}\times\mathbf{T}\right|^2}\\
    \end{split}
\end{equation}
and
\begin{equation}
\label{OmegaBOUNDARY}
    \begin{split}
        \mathit{\Omega}=2\pi -\oint\frac{\left|\mathbf{v}\right|\left(\mathbf{T}_\theta\cdot\left(\mathbf{v}\times\mathbf{T}\right)\right)}{\left|\mathbf{v}\times\mathbf{T}\right|^2}d\theta\mathrm{.}\\
    \end{split}
\end{equation}

The integral \eqref{OmegaBOUNDARY} is impossible to evaluate analytically for even the simplest of geometries, and so we approximate the boundary by a series of short straight lines. Because the calculations are not complicated a large number of boundary lines can be used to approximate smooth curves.

If the boundary in 3 dimensional space consists of straight lines which meet at `kinks' where their directions change, then the straight lines map to geodesics ($k_{\mathrm{g}}=0$) and
\begin{equation}
\label{angleSum}
\mathit{\Omega}  = 2\pi  - \sum\limits_i {{\alpha _i}}
\end{equation}
where $\alpha _i$ is the angle between two lines as seen on the surface of the unit sphere. This is equal to the angle between the two planes defined by the centre of the unit sphere and the each of the two straight line elements. Because the integration is very difficult for any curve beyond points on the axis of a circle, a curve is approximated by a series of straight lines and $\mathit{\Omega}$  is calculated numerically using the summation \eqref{angleSum}.
\begin{figure}[htp]
\centering
\includegraphics[width=0.5\textwidth]{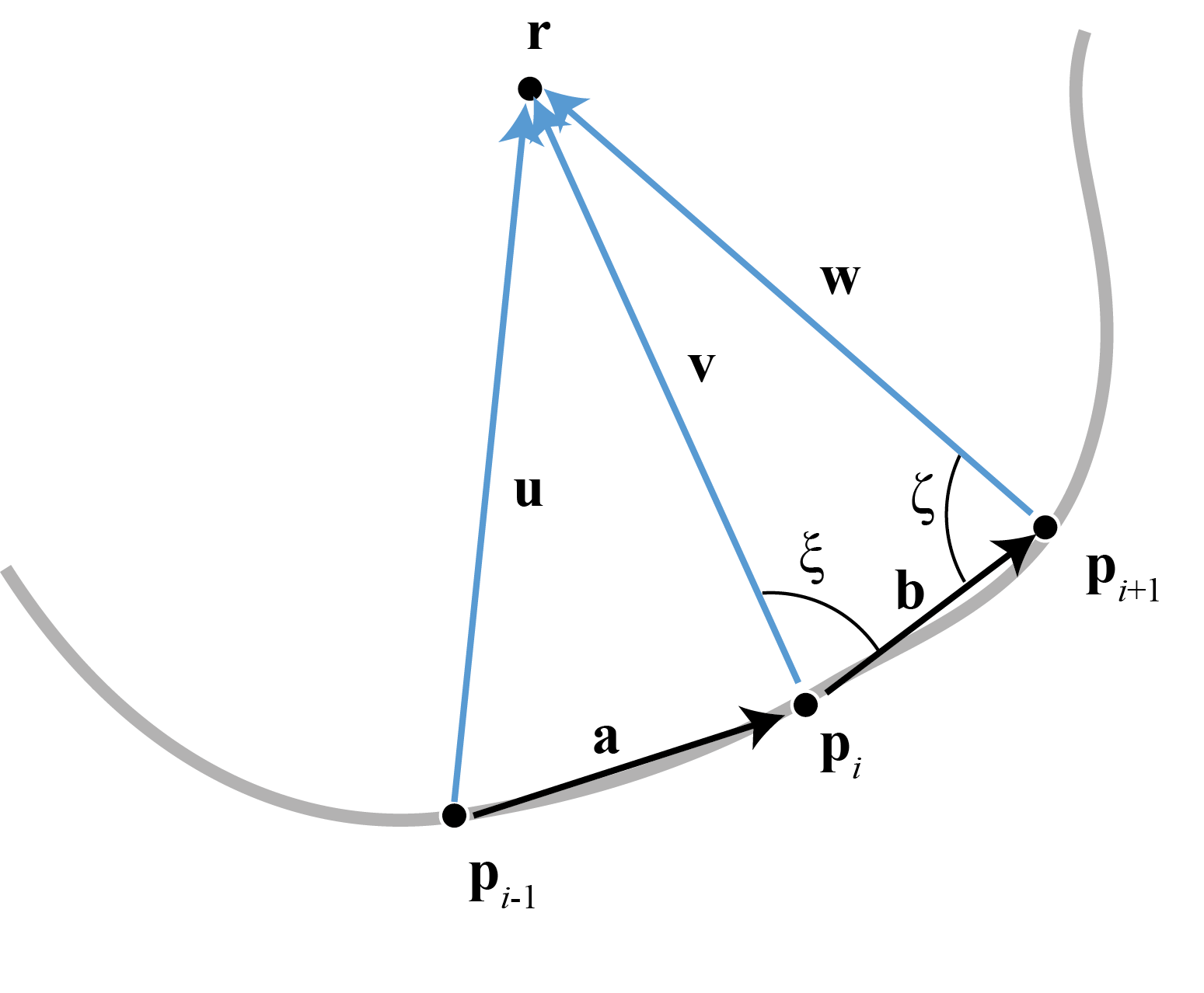}
\caption{Three points along the boundary and the point ${\mathbf{r}}$ at which to calculate the solid angle.}
\label{fig:SolidAngleContribution}
\end{figure}
Writing
\begin{align}
{\mathbf{a}} &= {{\mathbf{p}}_i} - {{\mathbf{p}}_{i - 1}}\\
{\mathbf{b}} &= {{\mathbf{p}}_{i + 1}} - {{\mathbf{p}}_i}\\
{\mathbf{v}} &= {\mathbf{r}} - {{\mathbf{p}}_i}
\end{align}
in which ${{\mathbf{p}}_{i - 1}}$, ${{\mathbf{p}}_{i}}$ and ${{\mathbf{p}}_{i + 1}}$ are the position vectors of 3 consecutive points on the boundary joined by straight lines and ${\mathbf{r}}$ is the point in space or apex at which to calculate the solid angle (Fig. \ref{fig:SolidAngleContribution}),
 $\alpha _i$ is calculated using
\begin{equation}
\sin {\alpha _i} = \frac{{\left[ {\left( {{\mathbf{a}} \times {\mathbf{v}}} \right) \times \left( {{\mathbf{b}} \times {\mathbf{v}}} \right)} \right] \cdot {\mathbf{v}}}}{{\left| {{\mathbf{a}} \times {\mathbf{v}}} \right|\left| {{\mathbf{b}} \times {\mathbf{v}}} \right|\left| {\mathbf{v}} \right|}}\mathrm{.}
\end{equation}
This follows from the fact that $\dfrac{{\mathbf{a}} \times {\mathbf{v}}}{\left| {{\mathbf{a}} \times {\mathbf{v}}} \right|}$ is a unit vector normal to the plane containing ${\mathbf{p}}_{i-1}$, ${\mathbf{p}}_i$ and ${\mathbf{r}}$, and $\dfrac{{\mathbf{b}} \times {\mathbf{v}}}{\left| {{\mathbf{b}} \times {\mathbf{v}}} \right|}$ is a unit vector normal to the plane containing ${\mathbf{p}}_i$, ${\mathbf{p}}_{i+1}$ and ${\mathbf{r}}$. Therefore $\dfrac{{\left( {{\mathbf{a}} \times {\mathbf{v}}} \right) \times \left( {{\mathbf{b}} \times {\mathbf{v}}} \right)}}{{\left| {{\mathbf{a}} \times {\mathbf{v}}} \right|\left| {{\mathbf{b}} \times {\mathbf{v}}} \right|}}$ is a vector in the direction of ${\mathbf{v}}$ whose magnitude is the sine of the angle between the two planes.

Using the properties of the vector triple product and the vector scalar product,
\begin{align}
\begin{split}
\sin {\alpha _i} &= \frac{{\left[ {{\mathbf{v}} \times \left( {{\mathbf{a}} \times {\mathbf{v}}} \right)} \right] \cdot \left( {{\mathbf{b}} \times {\mathbf{v}}} \right)}}{{\left| {{\mathbf{a}} \times {\mathbf{v}}} \right|\left| {{\mathbf{b}} \times {\mathbf{v}}} \right|\left| {\mathbf{v}} \right|}}
 = \frac{{\left[ {\left( {{\mathbf{v}} \cdot {\mathbf{v}}} \right){\mathbf{a}} - \left( {{\mathbf{a}} \cdot {\mathbf{v}}} \right){\mathbf{v}}} \right] \cdot \left( {{\mathbf{b}} \times {\mathbf{v}}} \right)}}{{\left| {{\mathbf{a}} \times {\mathbf{v}}} \right|\left| {{\mathbf{b}} \times {\mathbf{v}}} \right|\left| {\mathbf{v}} \right|}}\\
 &= \frac{{\left| {\mathbf{v}} \right|\left( {{\mathbf{a}} \times {\mathbf{b}}} \right) \cdot {\mathbf{v}}}}{{\left| {{\mathbf{a}} \times {\mathbf{v}}} \right|\left| {{\mathbf{b}} \times {\mathbf{v}}} \right|}}\mathrm{.}
\end{split}
\end{align}
It also gives
\begin{equation}
\cos {\alpha _i} = \frac{{\left( {{\mathbf{a}} \times {\mathbf{v}}} \right) \cdot \left( {{\mathbf{b}} \times {\mathbf{v}}} \right)}}{{\left| {{\mathbf{a}} \times {\mathbf{v}}} \right|\left| {{\mathbf{b}} \times {\mathbf{v}}} \right|}} = \frac{{\left( {{\mathbf{a}} \cdot {\mathbf{b}}} \right)\left( {{\mathbf{v}} \cdot {\mathbf{v}}} \right)} - {\left( {{\mathbf{a}} \cdot {\mathbf{v}}} \right)\left( {{\mathbf{b}} \cdot {\mathbf{v}}} \right)}}{{\left| {{\mathbf{a}} \times {\mathbf{v}}} \right|\left| {{\mathbf{b}} \times {\mathbf{v}}} \right|}}
\end{equation}
so that
\begin{equation}
\tan {\alpha _i} = \frac{{\left| {\mathbf{v}} \right|\left( {{\mathbf{a}} \times {\mathbf{b}}} \right) \cdot {\mathbf{v}}}}{{\left( {{\mathbf{a}} \cdot {\mathbf{b}}} \right)\left( {{\mathbf{v}} \cdot {\mathbf{v}}} \right)} - {\left( {{\mathbf{a}} \cdot {\mathbf{v}}} \right)\left( {{\mathbf{b}} \cdot {\mathbf{v}}} \right)}}\mathrm{.}
\label{tangent}
\end{equation}
\eqref{tangent} gives a positive value for $\alpha _i$ if $\mathbf{a}$ and $\mathbf{b}$ lie in the plane of the figure in Fig. \ref{fig:SolidAngleContribution} and if the point $\mathbf{r}$ lies above the plane of the figure as we look at it. Since $\alpha _i$ may lie between $-\pi$ and $\pi$, it is necessary to use the $\atantwo$ function in computer programming,

\begin{equation}
\begin{split}
    \alpha_i&= \atantwo\left(y,x\right)\\
    x&={\left( {{\mathbf{a}} \cdot {\mathbf{b}}} \right)\left( {{\mathbf{v}} \cdot {\mathbf{v}}} \right)} - {\left( {{\mathbf{a}} \cdot {\mathbf{v}}} \right)\left( {{\mathbf{b}} \cdot {\mathbf{v}}} \right)}\\
y&=\left| {\mathbf{v}} \right|\left( {{\mathbf{a}} \times {\mathbf{b}}} \right) \cdot {\mathbf{v}}\mathrm{.}\\
\end{split}
\label{atan2}
\end{equation}

\section{The special case of a triangle} \label{sec:SpecialCaseofTriangle}
In the case of a boundary consisting of a triangle with straight sides, application of \eqref{angleSum} and \eqref{tangent} gives \citep{10.2307/2691141}
\begin{equation}
\label{triangle}
\tan \frac{\mathit{\Omega} }{2} = \frac{{\left( {{\mathbf{u}} \times {\mathbf{v}}} \right) \cdot {\mathbf{w}}}}{{\left| {\mathbf{u}} \right|\left| {\mathbf{v}} \right|\left| {\mathbf{w}} \right| + \left( {{\mathbf{u}} \cdot {\mathbf{v}}} \right)\left| {\mathbf{w}} \right| + \left( {{\mathbf{v}} \cdot {\mathbf{w}}} \right)\left| {\mathbf{u}} \right| + \left( {{\mathbf{w}} \cdot {\mathbf{u}}} \right)\left| {\mathbf{v}} \right|}}\mathrm{.}
\end{equation}
where there are only 3 boundary points and ${\mathbf{u}}$, ${\mathbf{v}}$ and ${\mathbf{w}}$ are as shown in figure \ref{fig:SolidAngleContribution}.

In fact to get from  \eqref{angleSum} and \eqref{tangent} to \eqref{triangle} is not a trivial task, but \eqref{triangle} can be proven by observing that cutting a triangle into two parts through an apex produces two triangles and demonstrating that the solid angle subtended by the full triangle is the sum of the solid angle subtended by the 2 parts.

Writing
\begin{align}
\begin{split}
{\mathbf{p}}_1 &= L\left( {\cos \beta {\mathbf{i}} + \sin \beta {\mathbf{j}}} \right)\\
{\mathbf{p}}_2 &= 0\\
{\mathbf{p}}_3 &= L\left( {\cos \beta {\mathbf{i}} - \sin \beta {\mathbf{j}}} \right)\\
{\mathbf{r}} &= x{\mathbf{i}} + y{\mathbf{j}} + z{\mathbf{k}}
\end{split}
\end{align}
and letting $L\to \infty$, so that 2 of the corners of the triangle become infinitely far away, then
\begin{align}
\begin{split}
{\mathbf{u}} &=  - L\left( {\cos \beta {\mathbf{i}} + \sin \beta {\mathbf{j}}} \right)\\
{\mathbf{v}} &=  x{\mathbf{i}} + y{\mathbf{j}} + z{\mathbf{k}}\\
{\mathbf{w}} &= - L\left( {\cos \beta {\mathbf{i}} - \sin \beta {\mathbf{j}}} \right)
\end{split}
\end{align}
so that
\begin{equation*}
\begin{split}
    &\tan \frac{\mathit{\Omega} }{2}=\\
    &\frac{{2\cos \beta \sin \beta z}}{{\sqrt {{x^2} + {y^2} + {z^2}} \left( {1 + {{\cos }^2}\beta  - {{\sin }^2}\beta } \right) - \left( {x\cos \beta  - y\sin \beta } \right) - \left( {x\cos \beta  + y\sin \beta } \right)}}
\end{split}
\end{equation*}
which produces
\begin{equation}
\label{cone}
{x^2} + {y^2} + {z^2} = \frac{1}{{{{\cos }^2}\beta }}{\left( {\frac{{\sin \beta z}}{{\tan \dfrac{\mathit{\Omega} }{2}}} + x} \right)^2}\mathrm{.}
\end{equation}
This is the equation of a cone with an elliptic cross-section whose shape depends on the constants $\mathit{\Omega}$ and $\beta$.

Thus for any boundary shape, in the immediate region of any sharp corner the surface is locally equivalent to a cone with an elliptic cross-section. Of course, without the singularity of curvature at the apex of a cone, a smooth surface must lie in the plane defined by 2 straight lines where they meet.

Allowing $\beta \to 0$ in \eqref{cone} we obtain the special case of two parallel lines in which a circular cylinder replaces the elliptic cone, as we would expect from the inscribed angle theorem for a circle.
\section{The gradient of the solid angle subtended by a closed curve and the Biot-Savart law}
\label{solidAngleGradient}
We now return again to the case of a boundary of arbitrary shape. $\nabla\mathit{\Omega}$ is the normal to a surface of constant $\mathit{\Omega}$, see Fig. \ref{fig:gradientSolidAngle}. To find $\nabla\mathit{\Omega}$ it is necessary to imagine that the point in space ${\mathbf{r}}$ moves slightly and all the boundary points remain fixed.

From \eqref{Omega},
\begin{equation}
\nabla \mathit{\Omega}  =  - \sum\limits_i {\nabla {\alpha _i}}
\end{equation}
and its is possible to calculate ${\nabla {\alpha _i}}$ from \eqref{tangent} using
\begin{equation}
\nabla {\mathbf{v}} = \nabla {\left({{\mathbf{r}} - {\mathbf{p}}_i}\right)} = \nabla {\mathbf{r}} = {\mathbf{I}}
\end{equation}
in which $\mathbf{I}$ is the unit tensor defined by
\begin{equation*}
{\mathbf{I}} \cdot {\mathbf{f}} = {\mathbf{f}} \cdot {\mathbf{I}} = {\mathbf{f}}
\end{equation*}
where $\mathbf{f}$ is any vector. Thus, for example,
\begin{equation*}
\nabla {\left( {{\mathbf{b}} \cdot {\mathbf{v}}} \right)} = \nabla {\left( {{\mathbf{b}} \cdot {\mathbf{r}}} \right)} = {\mathbf{b}}\mathrm{,}
\end{equation*}
in which the vector $\mathbf{b}$ is taken as constant and only $\mathbf{r}$ and $\mathbf{v}$ vary.

After a not inconsiderable amount of working it is possible to obtain
\begin{equation}
    \begin{split}
    \label{delOmega}
\nabla\mathit{\Omega}&=-\sum\limits_i\left(\mathbf{v}-\mathbf{w}\right)\cdot\left(\frac{\mathbf{v}}{\left|\mathbf{v}\right|}-\frac{\mathbf{w}}{\left|\mathbf{w}\right|}\right)
\frac{\left(\mathbf{v}\times\mathbf{w}\right)}{\left|\mathbf{v}\times\mathbf{w}\right|^2}\\
    \end{split}
\end{equation}
in which again ${\mathbf{v}} = {\mathbf{r}} - {{\mathbf{p}}_i}$ and ${\mathbf{w}} = {\mathbf{r}} - {{\mathbf{p}}_{i + 1}}$ as shown in Fig. \ref{fig:SolidAngleContribution}.

Now consider the vector ${\mathbf{R}}$ from a typical point on the straight line between ${{\mathbf{p}}}_i$ and  ${{\mathbf{p}}}_{i + 1}$ to ${\mathbf{r}}$. The minimum magnitude of ${\mathbf{R}}$ is the perpendicular distance,
\begin{equation*}
h=\frac{\left| {\mathbf{v}}\times{\mathbf{w}} \right|}{\left| {\mathbf{v}}-{\mathbf{w}} \right|}\mathrm{.}
\end{equation*}
The contribution to $\nabla \mathit{\Omega}$ in \eqref{delOmega} is given by
\begin{equation}
\label{BS}
\begin{split}
-\int^{{\mathbf{w}}}_{{\mathbf{R}}={\mathbf{v}}}\frac{{d{\mathbf{R}}\times{\mathbf{R}}}}{\left|{\mathbf{R}}\right|^3
}
&=-\frac{{\mathbf{v}}\times{\mathbf{w}}}{\left| {\mathbf{v}}\times{\mathbf{w}} \right|}\int^{\pi-\zeta}_{\theta=\xi}\frac{\sin^3\theta h{d\left(h\cot\theta\right)}}{h^3}\\
&=-\frac{\left| {\mathbf{v}}-{\mathbf{w}} \right|\left({\mathbf{v}}\times{\mathbf{w}}\right)}{\left| {\mathbf{v}}\times{\mathbf{w}} \right|^2}\left(\cos\xi+\cos\zeta\right)\\
\end{split}
\end{equation}
where $\xi$ and $\zeta$ are the angles shown in Fig. \ref{fig:SolidAngleContribution}. It is easy to see that the final line of \eqref{BS} is identical to the quantity being summed in \eqref{delOmega}.

This (apart from a multiplying constant) is the Biot-Savart law familiar from potential theory where it gives the magnetic field due to a current carrying wire in magnetostatics or the fluid velocity due to a vortex in irrotational incompressible flow of a fluid. Magnetostatics is the study of magnetic fields produced by steady currents in wires, in which our $\mathit{\Omega}$ would correspond to the scalar potential, when multiplied by some constant. Similarly in fluid mechanics our $\mathit{\Omega}$ corresponds to the velocity potential.

Note that the production of solid angle is centred on the kinks between straight line segments of boundary whereas the gradient of the solid angle is produced by each line segment separately.

For numerical evaluation we can observe that
\begin{equation}
    \begin{split}
\left(\mathbf{v}-\mathbf{w}\right)\cdot\left(\frac{\mathbf{v}}{\left|\mathbf{v}\right|}-\frac{\mathbf{w}}{\left|\mathbf{w}\right|}\right)
&=\frac{\left(\left|\mathbf{v}\right|+\left|\mathbf{w}\right|\right)}{\left|\mathbf{v}\right|\left|\mathbf{w}\right|}\left(\left|\mathbf{v}\right|\left|\mathbf{w}\right|-\mathbf{v}\cdot\mathbf{w}\right)\\
\left|\mathbf{v}\times\mathbf{w}\right|^2&=\left|\mathbf{v}\right|^2\left|\mathbf{w}\right|^2-\left(\mathbf{v}\cdot\mathbf{w}\right)^2
    \end{split}
\end{equation}
so that we can rewrite \eqref{delOmega} as
\begin{equation}
\label{delOmegaIMPROVED}
    \begin{split}
\nabla\mathit{\Omega}&=-\sum\limits_i\frac{\left(\left|\mathbf{v}\right|+\left|\mathbf{w}\right|\right)\left(\mathbf{v}\times\mathbf{w}\right)}{\left|\mathbf{v}\right|\left|\mathbf{w}\right|\left(\left|\mathbf{v}\right|\left|\mathbf{w}\right|+\mathbf{v}\cdot\mathbf{w}\right)}
    \end{split}
\end{equation}
which is more accurate when $\mathbf{v}$ and $\mathbf{w}$ are almost parallel, in which case both the numerator and denominator are very small in \eqref{delOmega}.

Equation \eqref{delOmega} can be rewritten other ways, such as
\begin{equation}
 \label{delOmegaALTERNATIVE}
    \begin{split}
\nabla\mathit{\Omega}&=-\sum\limits_i
\frac{\left(\mathbf{a}\cdot\mathbf{v}\right)\left(\mathbf{a}\times\mathbf{v}\right)}{\left|\mathbf{v}\right|\left|\mathbf{a}\times\mathbf{v}\right|^2}
-\sum\limits_i
\frac{\left(\mathbf{b}\cdot\mathbf{v}\right)\left(\mathbf{b}\times\mathbf{v}\right)}{\left|\mathbf{v}\right|\left|\mathbf{b}\times\mathbf{v}\right|^2}\mathrm{,}\\
    \end{split}
\end{equation}
but \eqref{delOmegaIMPROVED} is to be preferred for numerical work.
\section{Curvature of a surface of constant solid angle} \label{sec:curvature}
We will not consider the curvature of a surface of constant solid angle in detail, although, there are various practical reasons why we should want to establish the curvature of a surface, such as structural analysis of a shell structure or the construction of a principal curvature grid on a surface. What follows is an outline of two alternatives on how one would find the curvature, but this section can be left by those only interested in establishing the surface itself, and not its curvature.

The first method utilises that $\nabla \mathit{\Omega}$ is normal to any surface $\mathit{\Omega} = \mathrm{constant}$ and ${\dfrac{{\nabla \mathit{\Omega} }}{{\left| {\nabla \mathit{\Omega} } \right|}}}$ is the unit normal. To find the curvature of a surface we need to find how the unit normal changes direction as we move across the surface.

In order to find the gradient of $\nabla\mathit{\Omega}$, that is the second order tensor $\nabla\nabla\mathit{\Omega}$, let us use \eqref{delOmegaALTERNATIVE} and remember that as the point $\mathbf{r}$ moves only $\mathbf{v}$ changes, while $\mathbf{a}$ and $\mathbf{b}$ are constant. Thus
\begin{equation}
\label{nablanabla}
    \begin{split}
\nabla\nabla\mathit{\Omega}=&-\sum\limits_i
\left(\frac{\mathbf{a}}{\left|\mathbf{v}\right|}
-\frac{\mathbf{v}\left(\mathbf{a}\cdot\mathbf{v}\right)}{\left|\mathbf{v}\right|^3}\right)\frac{\left(\mathbf{a}\times\mathbf{v}\right)}{\left|\mathbf{a}\times\mathbf{v}\right|^2}\\
&+\sum\limits_i
\frac{\left(\mathbf{a}\cdot\mathbf{v}\right)}{\left|\mathbf{v}\right|}\left(
\frac{\left(\mathbf{I}\times\mathbf{a}\right)}{\left|\mathbf{a}\times\mathbf{v}\right|^2}
-\frac{2\left(\left(\mathbf{a}\cdot\mathbf{v}\right)\mathbf{a}-\left(\mathbf{a}\cdot\mathbf{a}\right)\mathbf{v}\right)\left(\mathbf{a}\times\mathbf{v}\right)}{\left|\mathbf{a}\times\mathbf{v}\right|^4}
\right)\\
&\mathrm{plus\;a\;similar\;expression\;replacing\;}\mathbf{a}\mathrm{\;by\;}\mathbf{b}\\
    \end{split}
\end{equation}
in which we have used
\begin{equation}
    \begin{split}
        \left(\mathbf{I}\times\mathbf{a}\right)\cdot\left(\mathbf{a}\times\mathbf{v}\right)
        &=\mathbf{I}\cdot\left(\mathbf{a}\times\left(\mathbf{a}\times\mathbf{v}\right)\right)
        =\left(\mathbf{a}\cdot\mathbf{v}\right)\mathbf{a}-\left(\mathbf{a}\cdot\mathbf{a}\right)\mathbf{v}\mathrm{.}
    \end{split}
\end{equation}
and
\begin{equation}
    \begin{split}
        \mathbf{I}\times\mathbf{a}&=\left(\mathbf{ii}+\mathbf{jj}+\mathbf{kk}\right)\times\left(a_x\mathbf{i}+a_y\mathbf{j}+a_z\mathbf{k}\right)\\
        &=-a_z\left(\mathbf{ij}-\mathbf{ji}\right)
        -a_x\left(\mathbf{jk}-\mathbf{kj}\right)
        -a_y\left(\mathbf{ki}-\mathbf{ik}\right)\mathrm{,}
    \end{split}
\end{equation}
which can be found using (A.5.4) in \cite{Rubin2000}.

Examination of \eqref{nablanabla} confirms that the trace of $\nabla\nabla\mathit{\Omega}$,
\begin{equation}
    \nabla\cdot\nabla\mathit{\Omega}=\nabla^2\mathit{\Omega}=0
\end{equation}
so that $\mathit{\Omega}$ obeys Laplace's equation, as we would expect from potential theory.

We know that the second order tensor $\nabla\nabla\mathit{\Omega}$ in \eqref{nablanabla} must be symmetric but the quantity to be summed on the right hand side is not symmetric. One would expect that upon doing the summation the result would be symmetric.

The symmetric second order tensor
\begin{equation}
\boldsymbol{\upbeta} =  - \left( {{\mathbf{I}} - \frac{{\nabla \mathit{\Omega} \nabla \mathit{\Omega} }}{{{{\left| {\nabla \mathit{\Omega} } \right|}^2}}}} \right) \cdot \frac{\nabla \nabla \mathit{\Omega}}{\left| {\nabla \mathit{\Omega} } \right|}  \cdot \left( {{\mathbf{I}} - \frac{{\nabla \mathit{\Omega} \nabla \mathit{\Omega} }}{{{{\left| {\nabla \mathit{\Omega} } \right|}^2}}}} \right)
\end{equation}
tells how the direction of the unit normal varies as we move on the surface, that is in a direction perpendicular to the unit normal. $\boldsymbol{\upbeta}$ has no component normal to the surface and it is therefore a surface tensor and its components are known as the \emph{coefficients of the second fundamental form} in differential geometry \citep{Struik1961,Eisenhart1947,GreenZerna68}.

The two principal values of $\boldsymbol{\upbeta}$ are the principal curvatures and the corresponding principal directions are the principal curvature directions. 
Finding the principal curvature directions is important in practice for cladding a surface with plane quadrilaterals \citep{pottmann-2007-dsad} in which the edges of the quadrilaterals follow the principal curvature directions.

\begin{figure}[htp]
    \centering
    \includegraphics[width=1.0\textwidth]{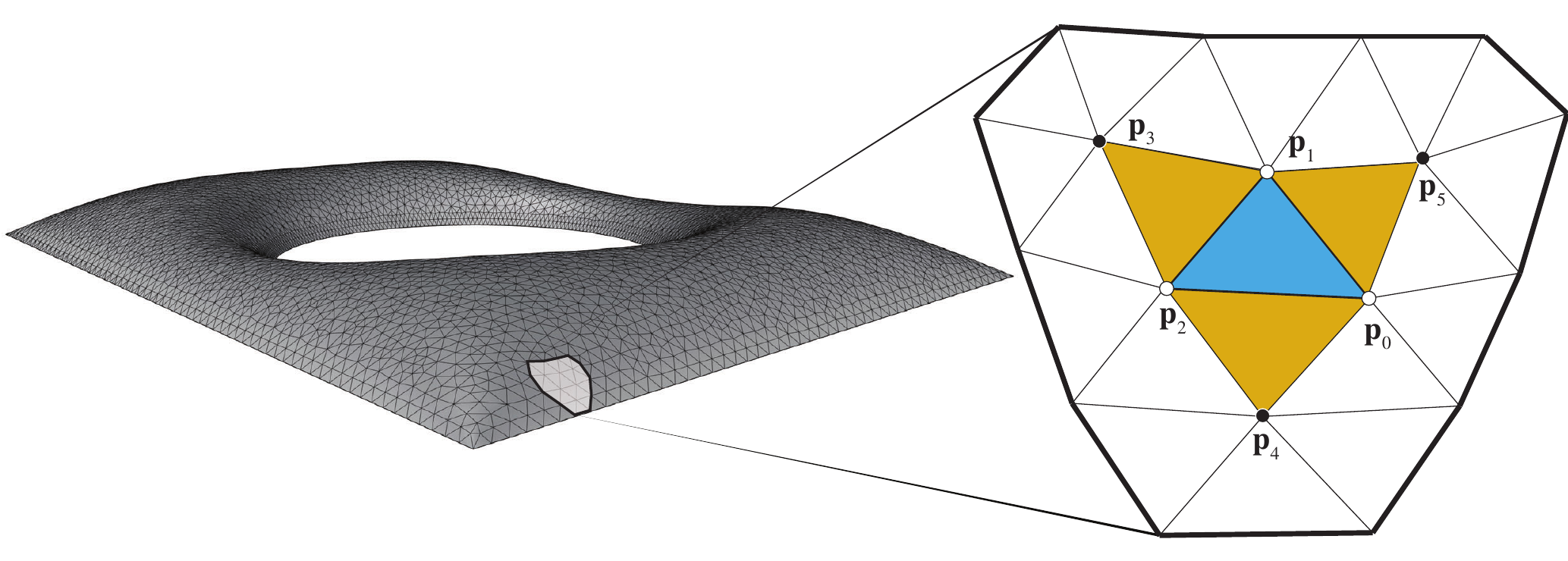}
    \caption{It is possible to get an approximate surface for the blue triangle passing trough the points $\mathbf{p}_i\left(x_i,y_i,z_i\right)$, $i=0$ to $5$ by solving the quadratic \eqref{eq:triangle}.}
    \label{fig:triangle}
\end{figure}

The second alternative, used to produce the principal curvature lines in the figures in section \ref{sec:examples}, requires a mesh in comparison to the previous method. For each mesh face it is possible to find a quadratic for $z$ in $x$ and $y$,
\begin{equation} \label{eq:triangle}
 z=ax^2+2bxy+cy^2+fx+gy+h,
\end{equation}
which passes through the six points $\mathbf{p}_i\left(x_i,y_i,z_i\right)$, $i=0$ to $5$, see Fig. \ref{fig:triangle}.
Having found the constants of the quadratic in \eqref{eq:triangle} it is straightforward to compute the principal curvature directions by solving the eigenvectors in \eqref{eq:principalDir2}.
\begin{equation}
\label{eq:principalDir2}
    \left[\left[\begin{array}{cc}
        \dfrac{\partial^2z}{\partial x^2} &  \dfrac{\partial^2z}{\partial x\partial y}\\
        \dfrac{\partial^2z}{\partial y\partial x} & \dfrac{\partial^2z}{\partial y^2}
    \end{array}\right]-\lambda\left[\begin{array}{cc}
        1+\left(\dfrac{\partial z}{\partial x}\right)^2 &  \dfrac{\partial z}{\partial x}\dfrac{\partial z}{\partial y}\\
        \dfrac{\partial z}{\partial y}\dfrac{\partial z}{\partial x} & 1+\left(\dfrac{\partial z}{\partial y}\right)^2
    \end{array}\right]\right]
    \left[\begin{array}{c}
         v_x  \\\\
         v_y 
    \end{array}\right]
\end{equation}

\section{Multiple boundaries}
The current in a single wire must have the same magnitude at all points, and the same applies to the strength of a single vortex. Thus with a single wire the solid angle subtended by the wire at any point is some constant times the scalar potential. If we have more than one wire then we have to include the value of the current since it is possible to have different electric currents, or their equivalents, in the different wires. It is also possible to apply Kirchhoff's current law where sections of boundary meet.

Then \eqref{Omega} becomes
\begin{equation}
\label{OmegaMULTIPLE}
\mathit{\Omega}=\sum_w I_w\left(2\pi-\int\limits_{\partial A_w}{{k_{\mathrm{g}}}ds}\right)
\end{equation}
in which the contributions of the different parts of the boundary are weighted by the currents $I_w$ in each of the wires $w$. Note that in applying Kirchhoff's current law we think of each wire $w$ being a closed loop and along certain lengths two or more wires may run alongside each other in which case the combined current is the sum of the individual current loops, which may be positive or negative.

Changing the sign of $I_w$ reverses the current, which should have exactly the same effect as reversing the direction of the integral $\int\limits_{\partial A_w}{{k_{\mathrm{g}}}ds}$. However, we then have to be careful to take into account what happens to the $2\pi$ in \eqref{OmegaMULTIPLE}. To get over this, it is often better to get rid of the $2\pi$ and instead use
\begin{equation}
\label{OmegaMULTIPLEno2pi}
\mathit{\Omega}=-\sum_w\left(I_w\int\limits_{\partial A_w}{{k_{\mathrm{g}}}ds}\right)
\end{equation}
to define our constant $\mathit{\Omega}$ surface. Even so we have to be aware that there is an uncertainty of $4\pi$ in each integral $\int\limits_{\partial A_w}{{k_{\mathrm{g}}}ds}$.
\section{Examples}
\label{sec:examples}
\begin{figure}[htp]
    \centering
    \includegraphics[width=0.8\textwidth]{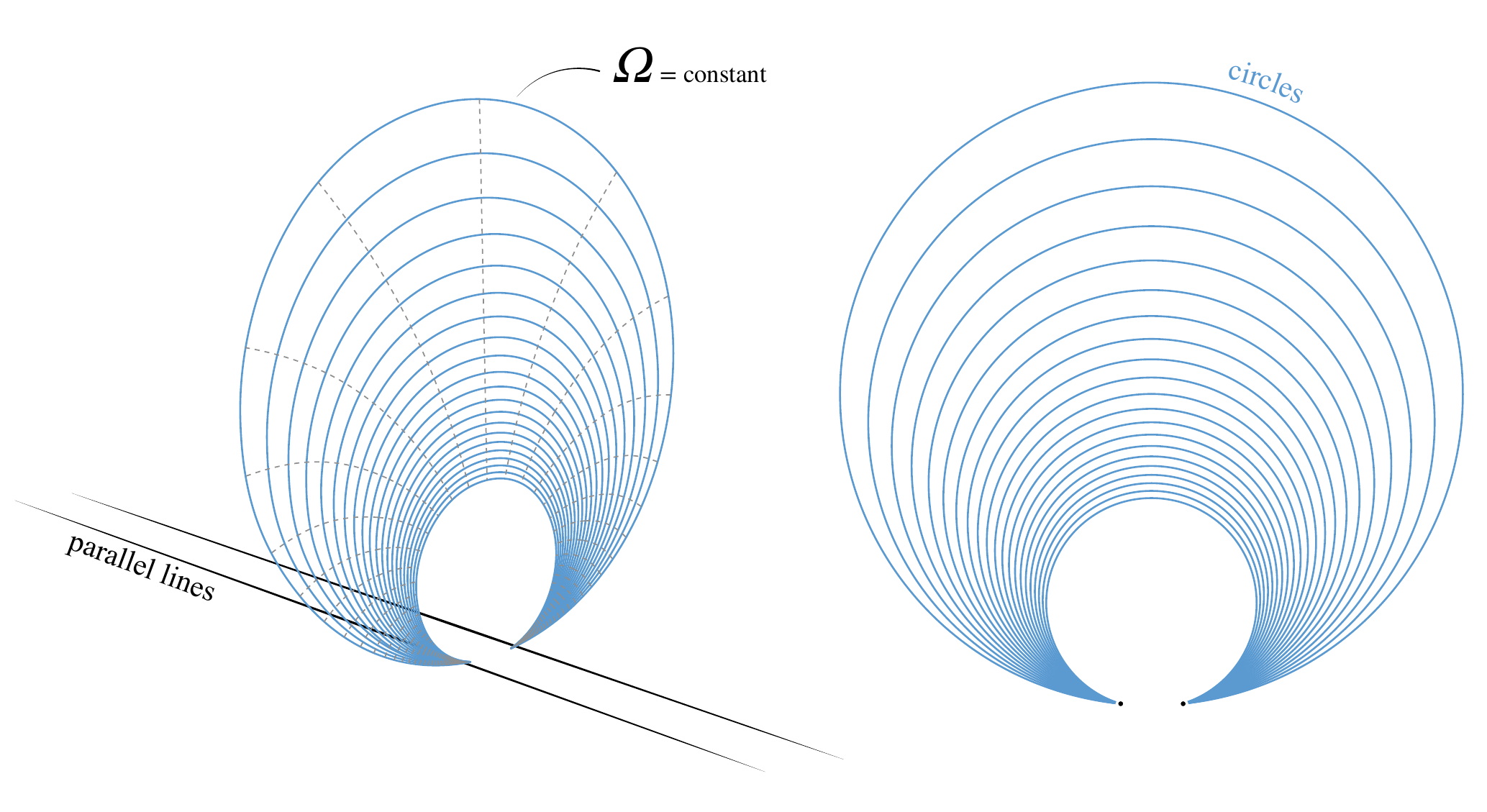}
    \caption{Potential flow produced by two straight parallel equal and opposite vortices modeled by a long rectangle using our method.}
    \label{fig:conformalMap}
\end{figure}

There are few analytical solutions of constant solid angle surfaces to compare with a numerical implementation. However, as described in Fig. \ref{sec:SpecialCaseofTriangle} a boundary consisting of two infinitely long parallel lines produces circular surfaces of constant solid angle. This is the equivalent to the lines of constant velocity potential produced by two infinite straight parallel equal and opposite line vortices \citep{Lamb1932}. We used a long thin rectangle to produce Fig. \ref{fig:conformalMap}, where the blue curves are indeed circles.

\begin{figure}[htp]
    \centering
    \includegraphics[width=0.8\textwidth]{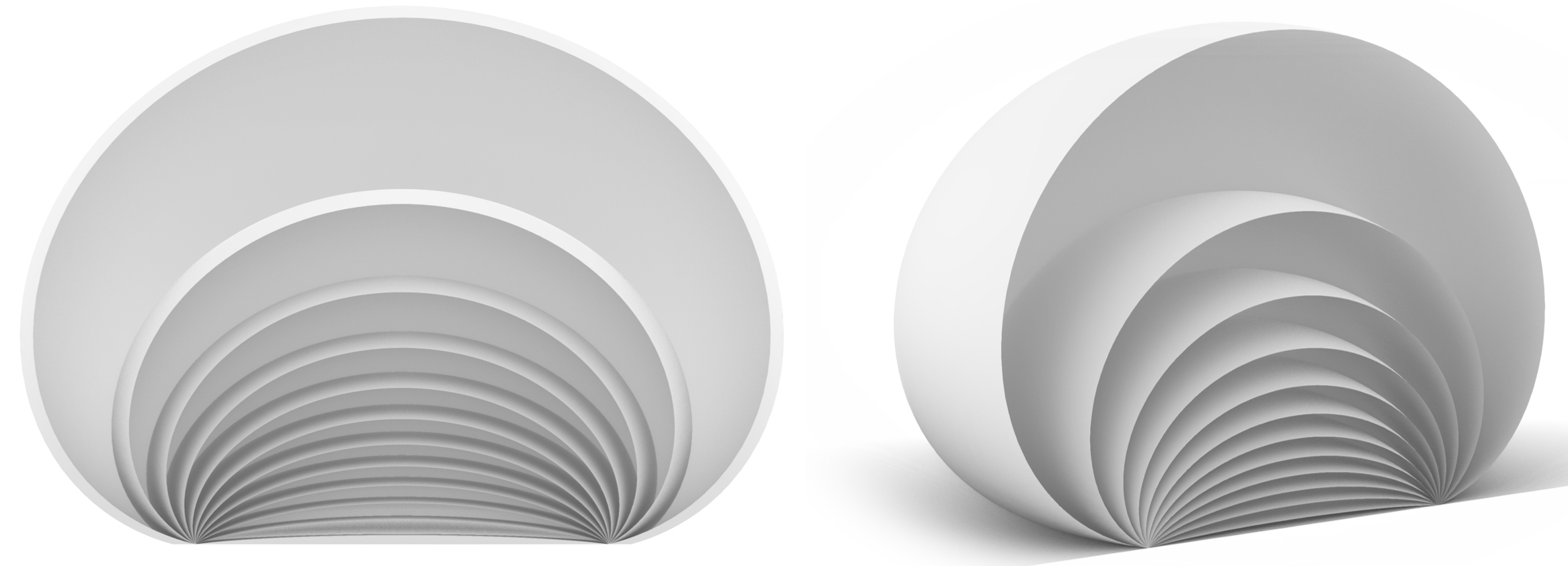}
    \caption{Surfaces with different constant solid angles on a circular boundary.}
    \label{fig:circularboundary_1}
\end{figure}

On a circular boundary, as in Fig. \ref{fig:circularboundary_1}, the shape is not a sphere since it is more flat on the top with increasing curvature near the boundary. The shape is more similar to a water droplet on a flat surface deformed by gravity than a sphere.

\begin{figure}[htp]
    \centering
    \includegraphics[width=0.7\textwidth]{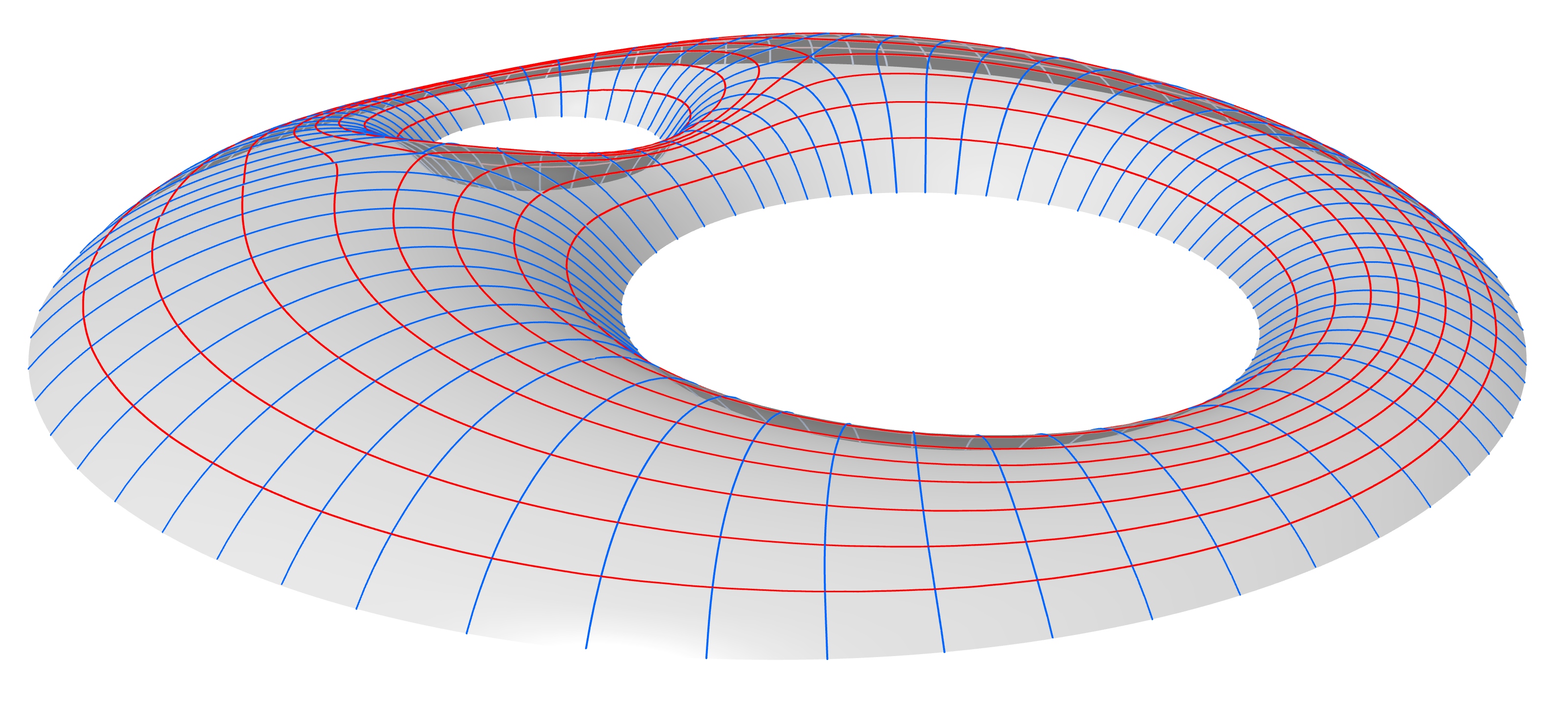}
    \caption{ Elevation of surface with constant solid angle having the the same boundary curves as Figs. \ref{fig:solidangleSurface} and \ref{fig:patternTop}, three circles with varying diameter that lies the same plane. The red and blue curves follows the principal curvature directions.}
      \label{fig:PatternSide}
\end{figure}

\begin{figure}[htp]
    \centering
    \includegraphics[width=1\textwidth]{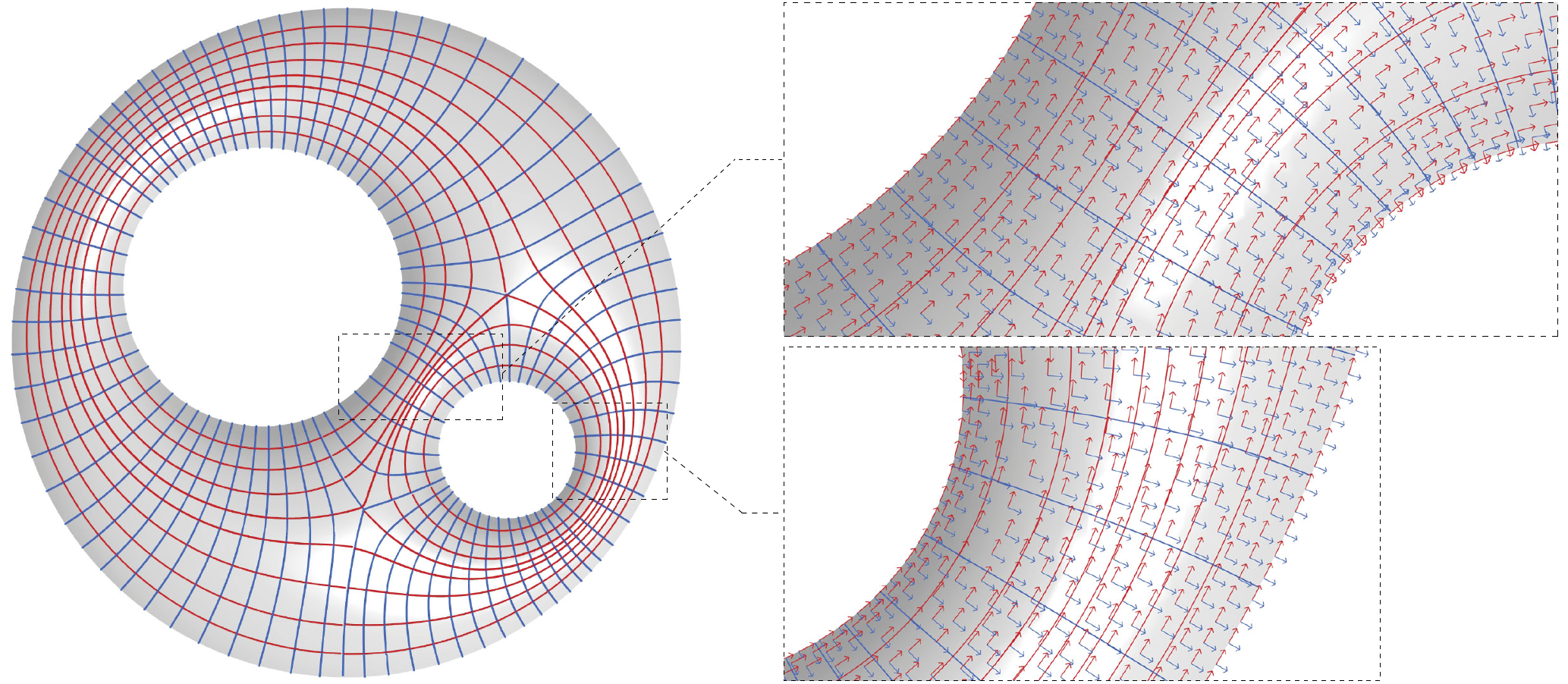}
    \caption{Top view of surface with constant solid angle having the the same boundary curves as Figs. \ref{fig:solidangleSurface} and \ref{fig:PatternSide}. Red and blue curves follows the principal curvature directions. Two areas are zoomed in to highlight that the principal curvature directions, seen as vectors, follow the surface boundary.}
    \label{fig:patternTop}
\end{figure}

Figures \ref{fig:PatternSide} and \ref{fig:patternTop} have the same boundary curves as the surface in Fig. \ref{fig:solidangleSurface}, three circles with different radius in the same plane, and the red and blue curves follows the principal curvature lines. In Fig. \ref{fig:patternTop} one can see in the zoomed in areas that the  principal curvature directions follow the surface boundary.   

\begin{figure}[htp]
    \centering
    \includegraphics[width=1.0\textwidth]{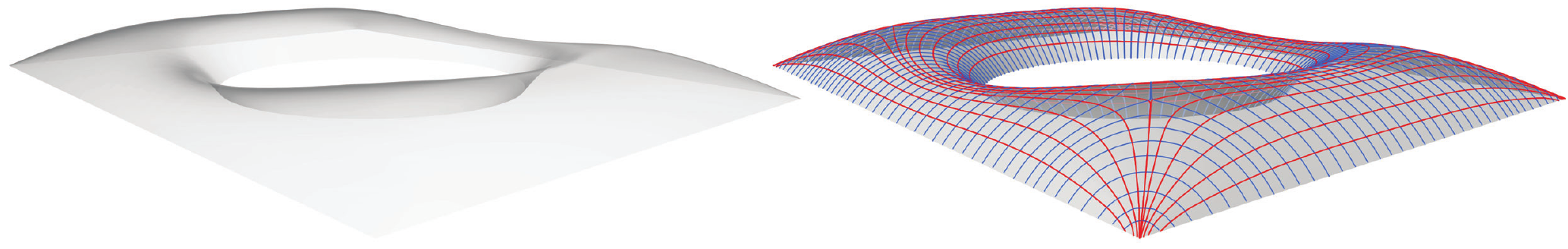}
    \caption{Elevation of surface with constant solid angle having the the same boundary curves as the British Museum Great Court roof. It is the same surfaces as in  seen in Figs. \ref{fig:britishPrincipalTop} and \ref{fig:britishPrincipalElavation} }
    \label{fig:britishPrincipalPersp}
\end{figure}

\begin{figure}[htp]
    \centering
    \includegraphics[width=1.0\textwidth]{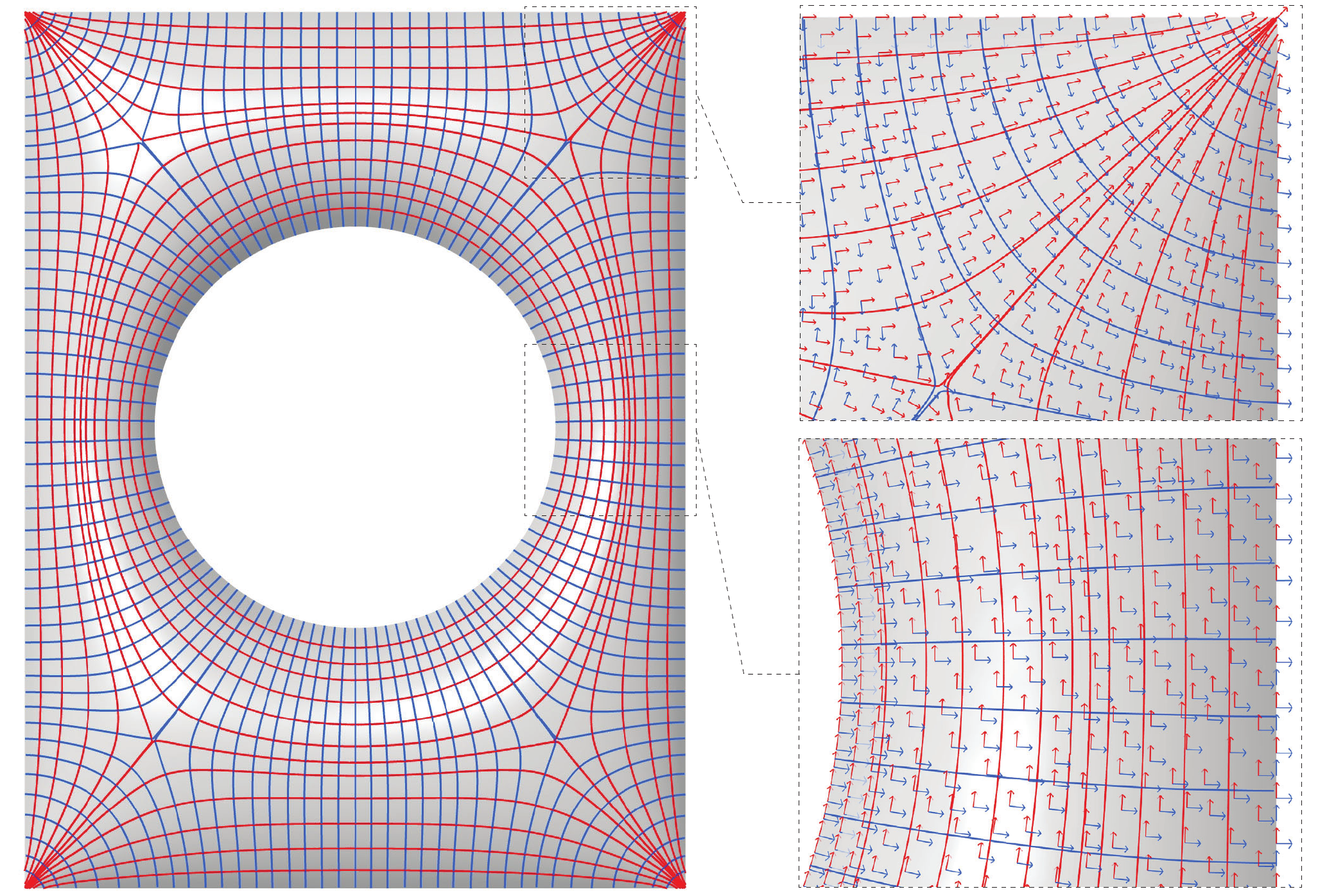}
    \caption{  Top view of the same surfaces seen in Figs. \ref{fig:britishPrincipalPersp} and \ref{fig:britishPrincipalElavation}. The red and blue curves are the principal curvature lines. Two areas are zoomed in to highlight that the principal curvature directions, seen as vectors, follow the surface boundary.}
    \label{fig:britishPrincipalTop}
\end{figure}

\begin{figure}[htp]
    \centering
    \includegraphics[width=1.0\textwidth]{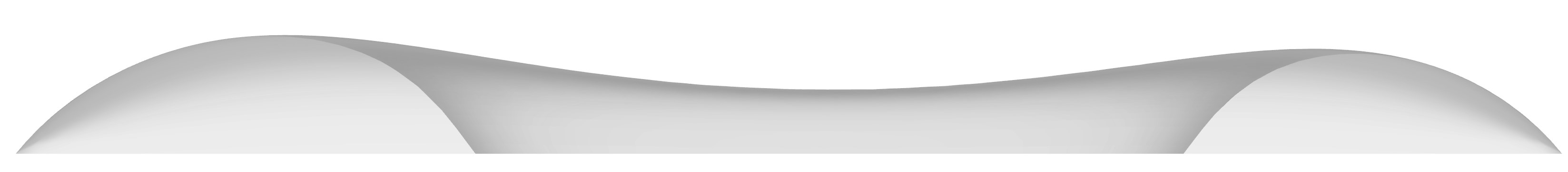}
    \caption{Elevation of surface with constant solid having the the same boundary curves as the British Museum Great Court roof. It is the same surfaces as in  seen in in Figs. \ref{fig:britishPrincipalPersp} and \ref{fig:britishPrincipalTop}.The material is transparent to visualise the interior space.} 
    \label{fig:britishPrincipalElavation}
\end{figure}

In Figs. \ref{fig:britishPrincipalPersp} to \ref{fig:britishPrincipalElavation} we have chosen the same boundary curves as the British Museum Great Court roof \cite{Williams2001} (Fig. \ref{fig:britishCourt}(a)) and placed them in the same plane. The principal curvature lines follow both the circular and the rectangular boundaries. Only in the close proximity of the corners can one see that the principal curvature directions diverge towards the corner. 


\begin{figure}[htp]
    \centering
    \includegraphics[width=0.8\textwidth]{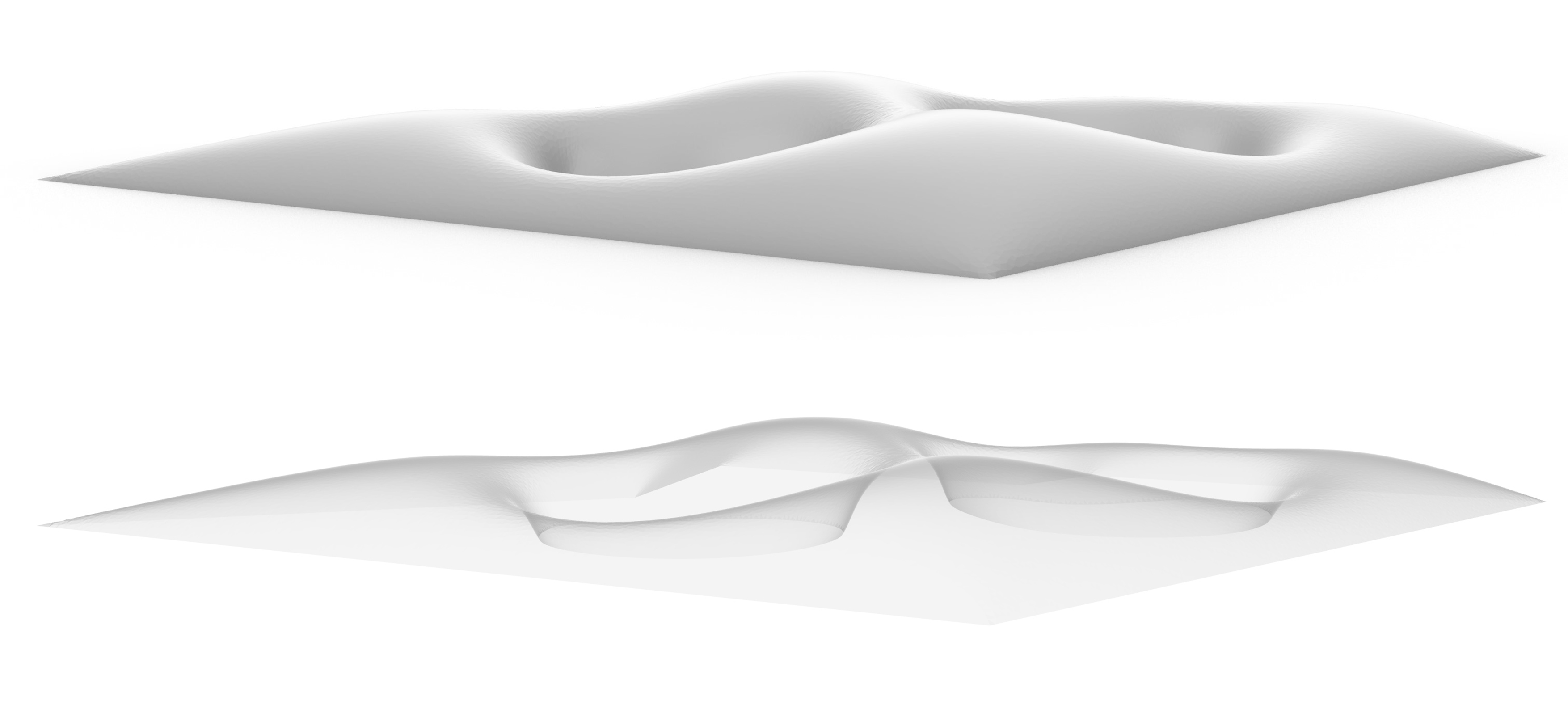}
    \caption{Perspective views of the constant solid angle surface in Figs. \ref{fig:complexsurftop} and \ref{fig:complexsurfelavation}. The slope is constant around each boundary, but varies from boundary to boundary.}
    
    \label{fig:complexsurfsec}
\end{figure}

\begin{figure}[htp]
    \centering
    \includegraphics[width=0.8\textwidth]{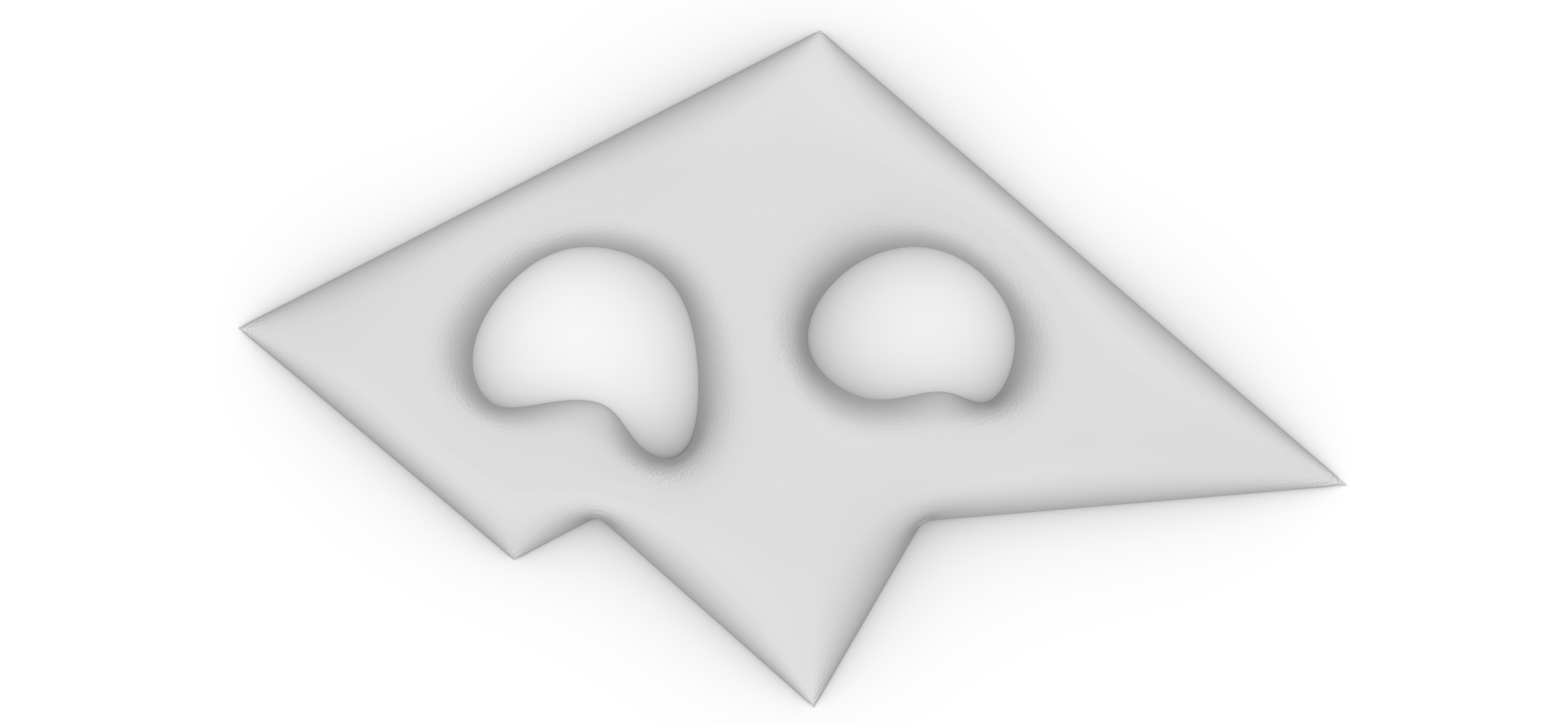}
    \caption{ Top view of the surface seen in Figs. \ref{fig:complexsurfsec} and \ref{fig:complexsurfelavation}. }
    \label{fig:complexsurftop}
\end{figure}

\begin{figure}[htp]
    \centering
    \includegraphics[width=1.0\textwidth]{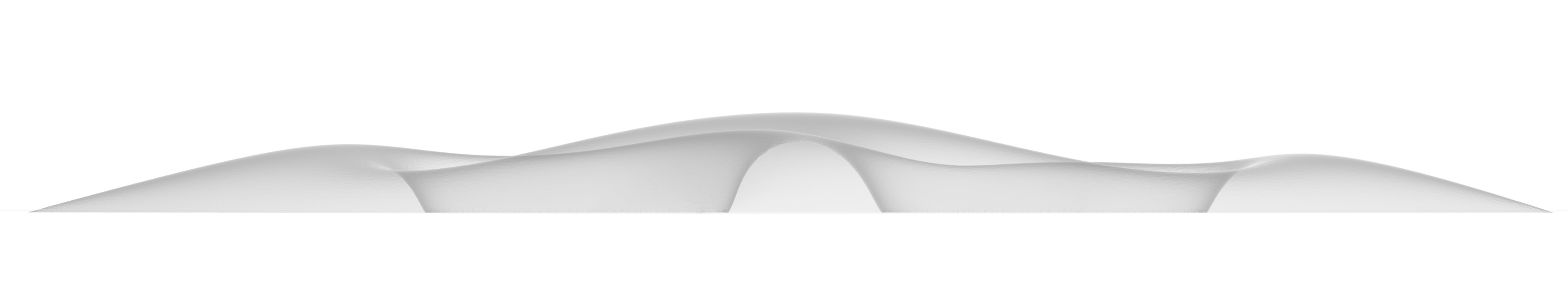}
    \caption{ Elevation of the surface in Figs. \ref{fig:complexsurfsec} and \ref{fig:complexsurftop}. The material is transparent to visualise the interior space}
    \label{fig:complexsurfelavation}
\end{figure}

Figures \ref{fig:complexsurfsec} to \ref{fig:complexsurfelavation} are all the same surface whose three boundary curves lie in the same plane. The boundary conditions were chosen to be quite problematic, but even so, the slope is constant along the boundary curves with a fine mesh.

\begin{figure}[htp]
    \centering
    \begin{tikzpicture}
    \node (myfirstpic) at (0,0)    
    {\includegraphics[width=0.95\textwidth]{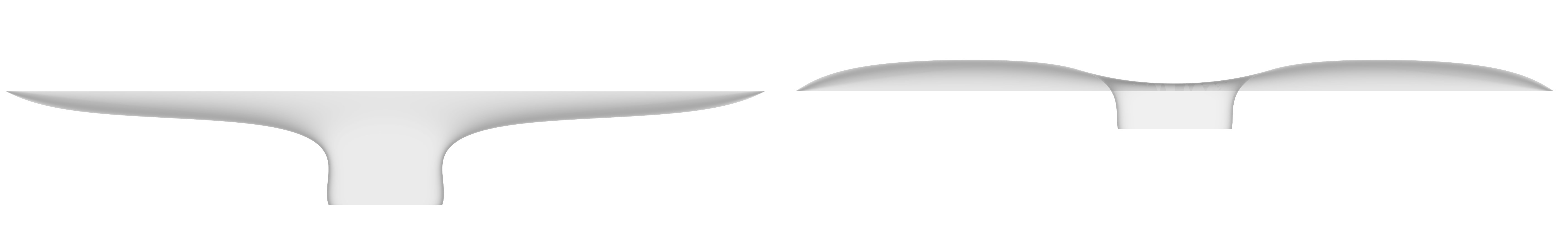}};
    \path (-7,-0.5) node[below] {(a)};
    \path (0.5,-0.5) node[below] {(b)};
    \end{tikzpicture}
    \caption{ Elevation of two surfaces with the same boundary curves, a rectangle and a circle which is lowered in relation to rectangle. In (a) the circle is placed lower than that in (b).}
    \label{fig:britishDifferentRatio}
\end{figure}

\begin{figure}[htp]
    \centering
     \begin{tikzpicture}
    \node (myfirstpic) at (0,0)    
    {\includegraphics[width=0.95\textwidth]{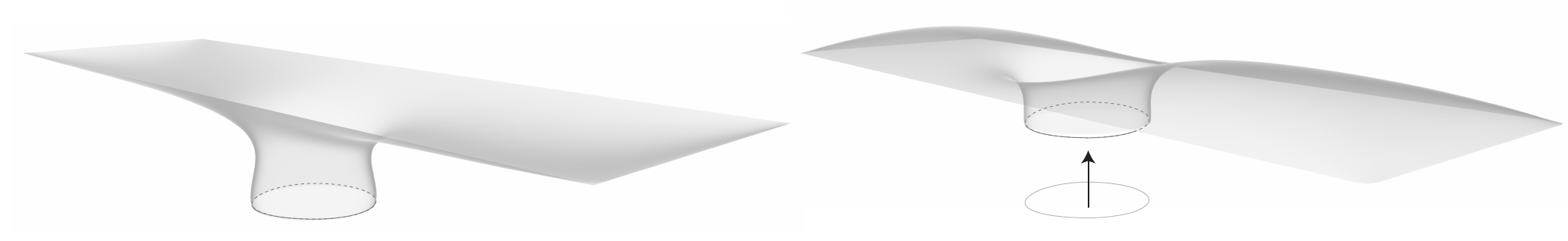}};
    \path (-7,-0.5) node[below] {(a)};
    \path (0.5,-0.5) node[below] {(b)};
    \end{tikzpicture}
    \caption{Perspective view of the same geometry as in Figs. \ref{fig:britishDifferentRatio} (a) and \ref{fig:britishDifferentRatio} (b)  rendered with a transparent material.}
    \label{fig:amancio}
\end{figure}

The surfaces in Figs. \ref{fig:britishDifferentRatio} and \ref{fig:amancio} have the same boundary curves, a rectangular exterior curve and a circular interior curve. Because the two boundary curves are at different levels, the slope around each curve will not be constant. Nevertheless we can rotate the slope around each curve independently by varying the current in the equivalent wire. 

\begin{figure}[htp]
    \centering
    \begin{tikzpicture}
    \node (myfirstpic) at (0,0)
    {\includegraphics[width=0.9\textwidth]{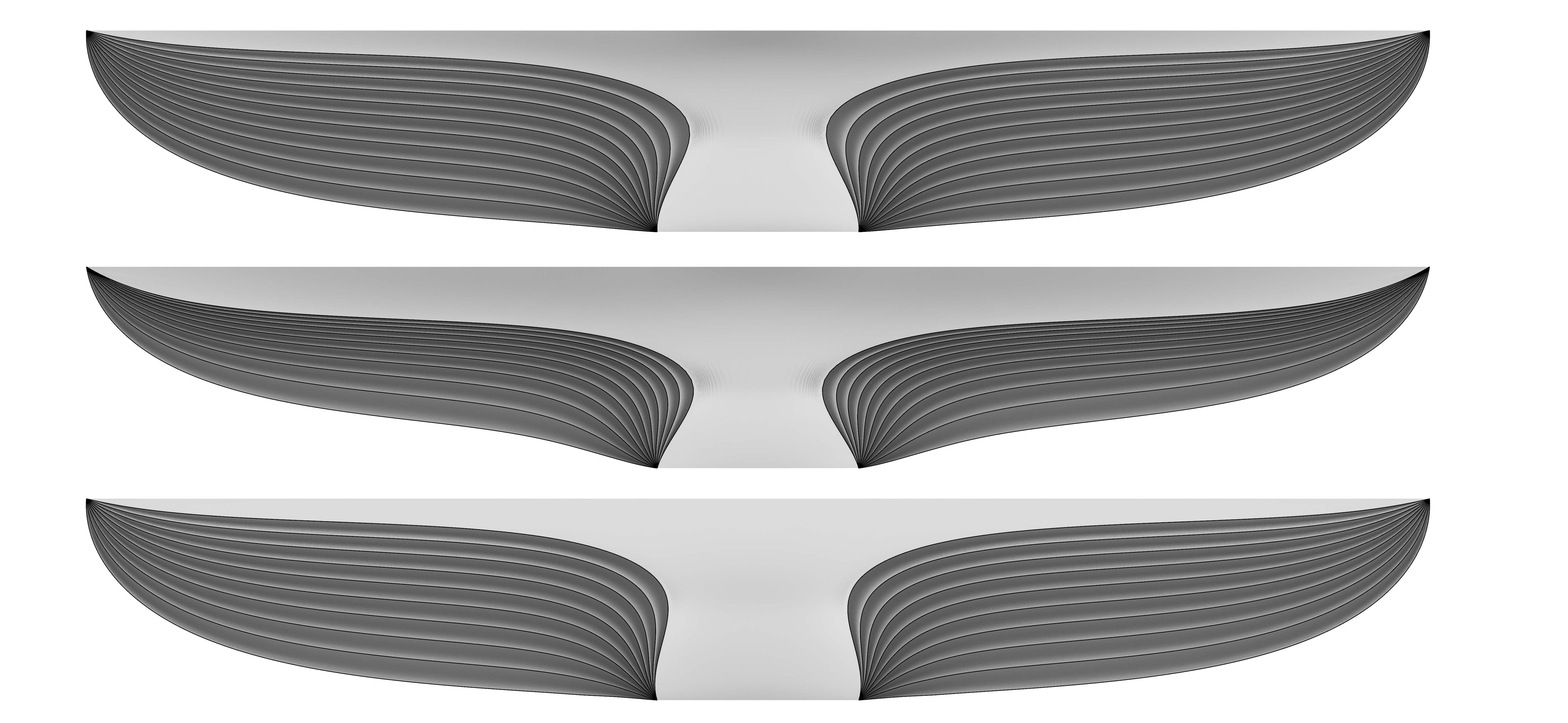}};
    \label{fig:amanciovarying}
    \path (-7,2.5) node[below] {(a)};
    \path (-7,0.5) node[below] {(b)};
    \path (-7,-1.5) node[below] {(c)};
    \end{tikzpicture}
    \caption{These three sections illustrate the effect of varying the solid angle and the ratio of the current in the circular and rectangular wires. The boundary conditions are the same as in Fig. \ref{fig:britishDifferentRatio} (a). In (a) the value of the solid angle varies while the current is the same in both wires. In (b) the current in the rectangular wire is varying for a value of constant solid angle, while in (c) the current in the circular wire is varying.}
\label{fig:amanciovaryingSECTIONS}
\end{figure}
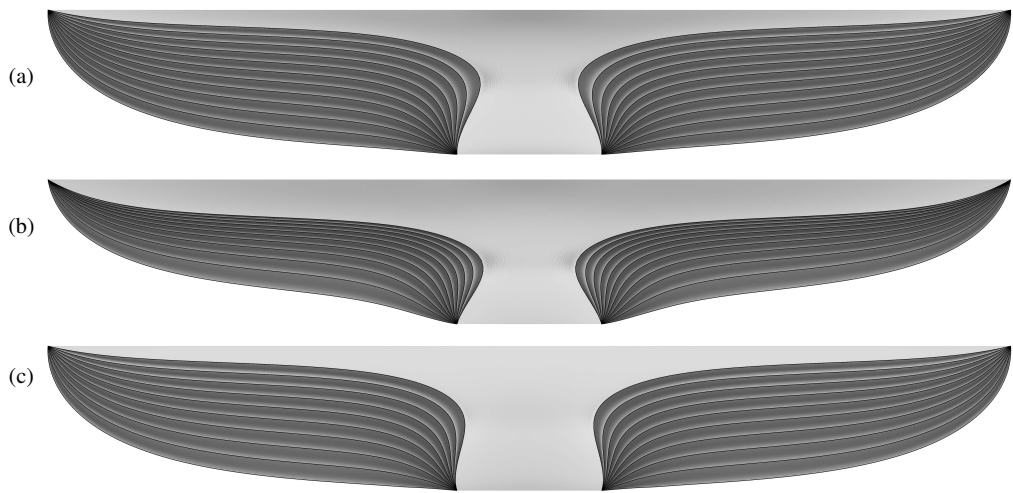

\begin{figure}[htp]
    \centering
     \begin{tikzpicture}
    \node (myfirstpic) at (0,0)
    {\includegraphics[width=1.0\textwidth]{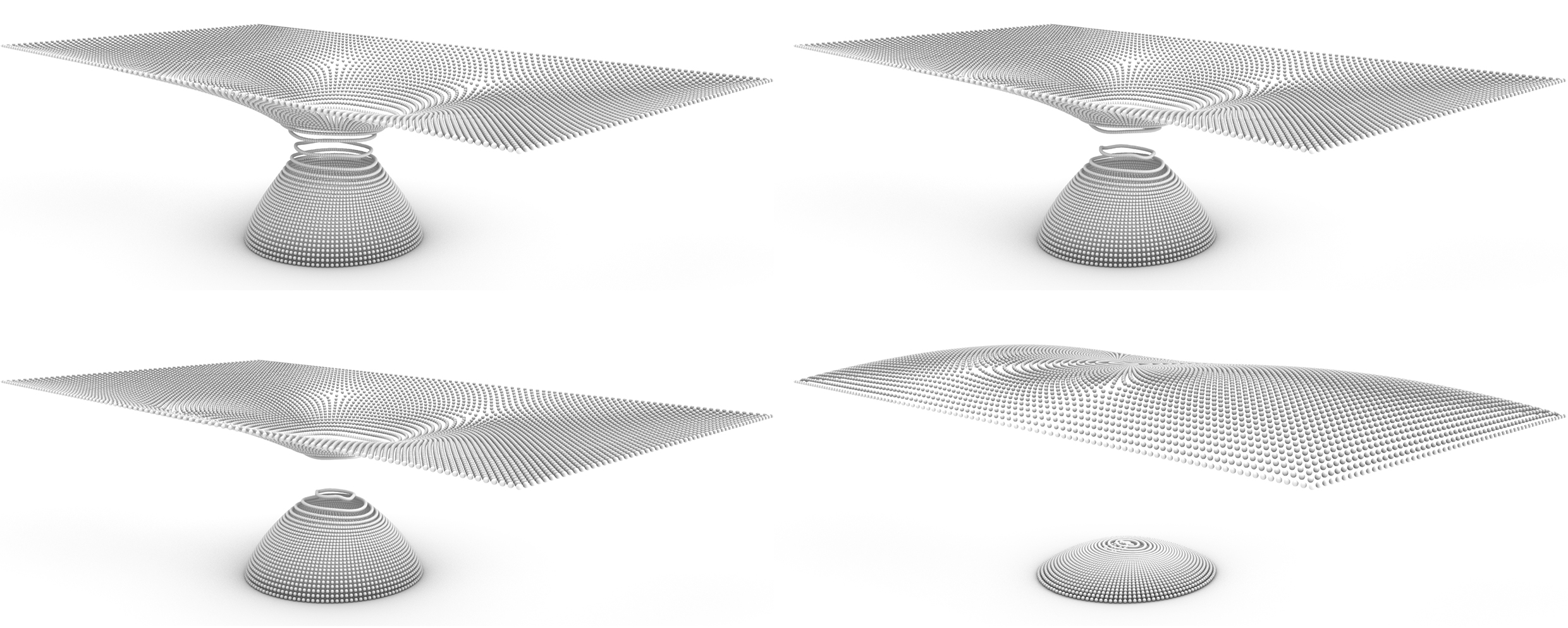}};
    \path (-7,1) node[below] {(a)};
    \path (-7,-2.25) node[below] {(c)};
    \path (1,1) node[below] {(b)};
    \path (1,-2.25) node[below] {(d)};
     \end{tikzpicture}
    \caption{Surfaces are represented by points with varying values of solid angle. The boundary curves are the similar to Fig. \ref{fig:britishDifferentRatio} (a), a rectangle and a circle in horizontal planes at different levels . At a certain value of the constant solid angle, the surface in (a) starts to separate into two surfaces seen in (b) and (c). In (d) the surface is split in two surfaces, both having the same constant solid angle.}
    \label{fig:surface_seperation}
\end{figure}
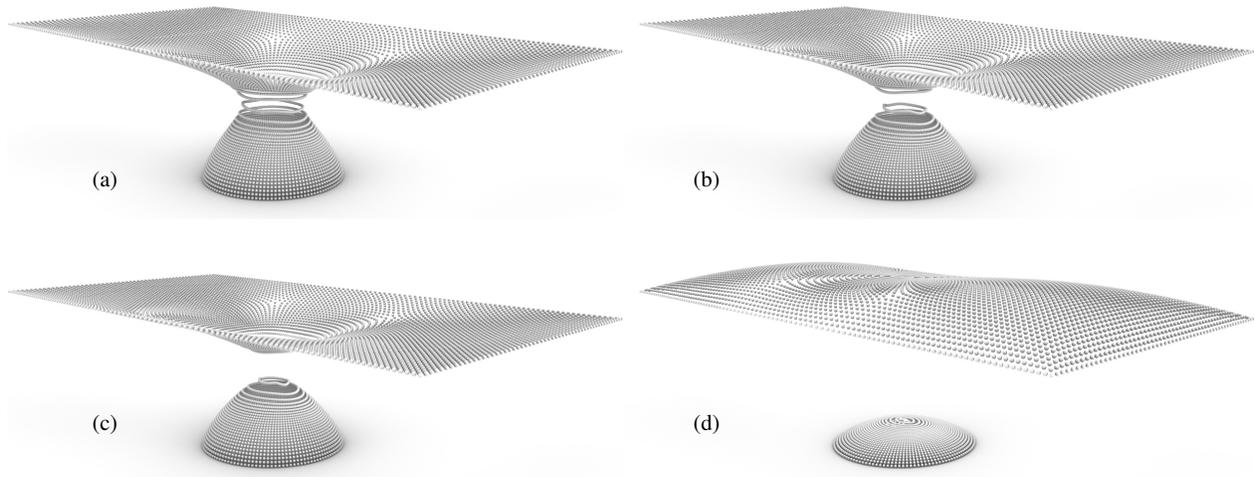

The relative height difference of the interior curve between Fig. \ref{fig:britishDifferentRatio} (a) and Fig. \ref{fig:britishDifferentRatio} (b) makes the slope at the exterior curve either positive or negative. Thus,  Fig. \ref{fig:britishDifferentRatio} (a) and Fig. \ref{fig:amancio} (a)  are reminiscent of the umbrella shells by Felix Candela and Amancio Williams \citep{Rian2014}, while Fig. \ref{fig:britishDifferentRatio} (b) and Fig. \ref{fig:amancio} (b) resemble the shape of the British Museum Great Court roof. Changing the value of the solid angle will not make the slope change from positive to negative as seen in Figs. \ref{fig:amanciovaryingSECTIONS} (a) to (c), unless breaking the surface into two. \newline \indent
The phenomenon of surface separation is illustrated in Fig. \ref{fig:surface_seperation} using points rather than a mesh and the same boundary conditions as Fig. \ref{fig:britishDifferentRatio} (a). One can see the surface as one in Fig. \ref{fig:surface_seperation} (a) but by changing the constant solid angle gradually one can see the surface starting to separate in Figs. \ref{fig:surface_seperation} (b) to (c). By further changing the solid angle the surface separate completely into two independent surfaces in Fig. \ref{fig:surface_seperation} (d), both having the same constant solid angle.

\begin{figure}[htp]
    \centering
    
    {\includegraphics[width=0.95\textwidth]{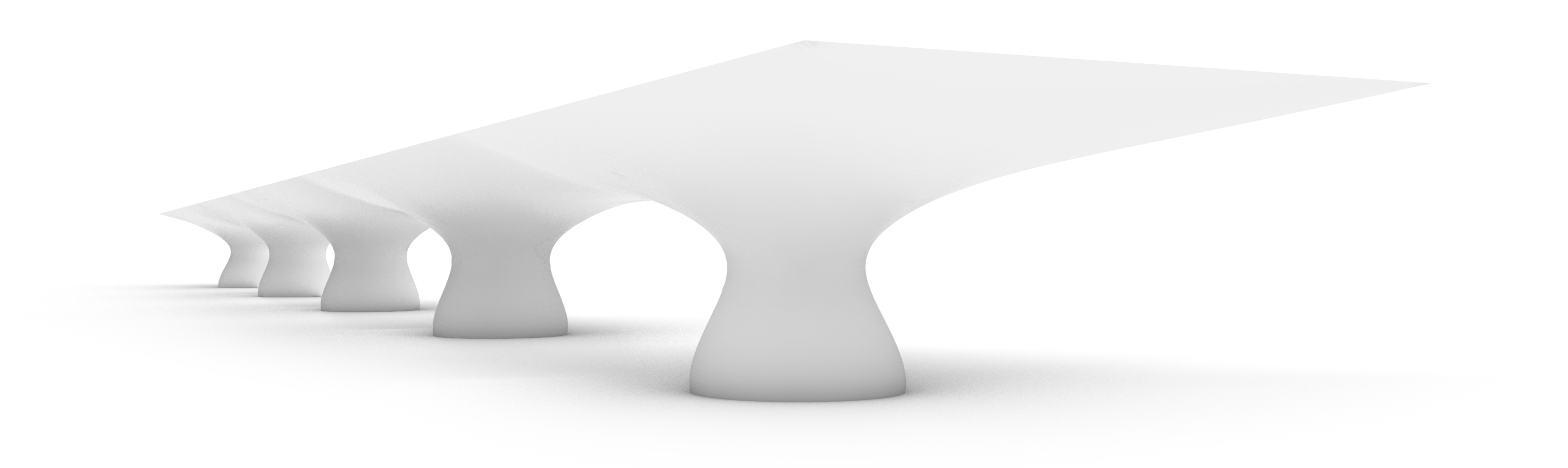}}
   
    \caption{ Perspective of a bridge design with a constant solid angle surface.}
    \label{fig:bridgeperspective}
\end{figure}

\begin{figure}[htp]
    \centering
    
    {\includegraphics[width=0.95\textwidth]{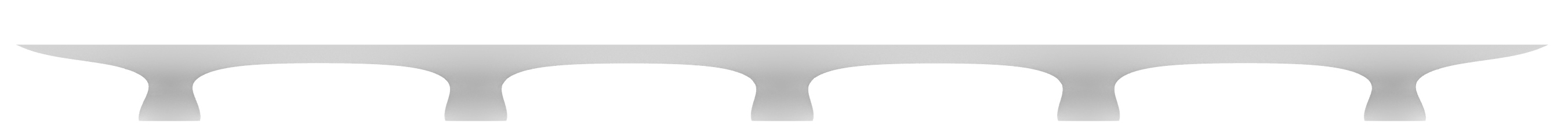}}
   
    \caption{ Elevation of a bridge design with a constant solid angle surface.}
    \label{fig:bridgeelavation}
\end{figure}

\begin{figure}[htp]
    \centering
    {\includegraphics[width=0.95\textwidth]{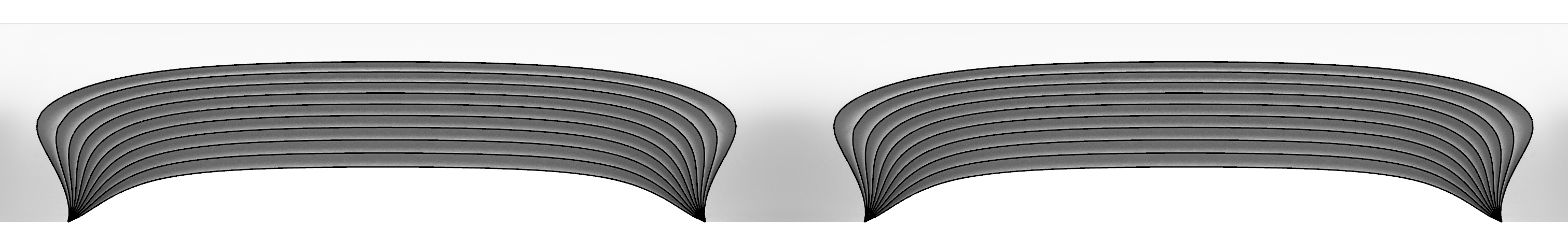}}
    \caption{Longitudinal sections of bridge designs with varying value of solid angle.}
    \label{fig:bridgsection}
\end{figure}

\begin{figure}[htp]
    \centering
    \includegraphics[width=0.8\textwidth]{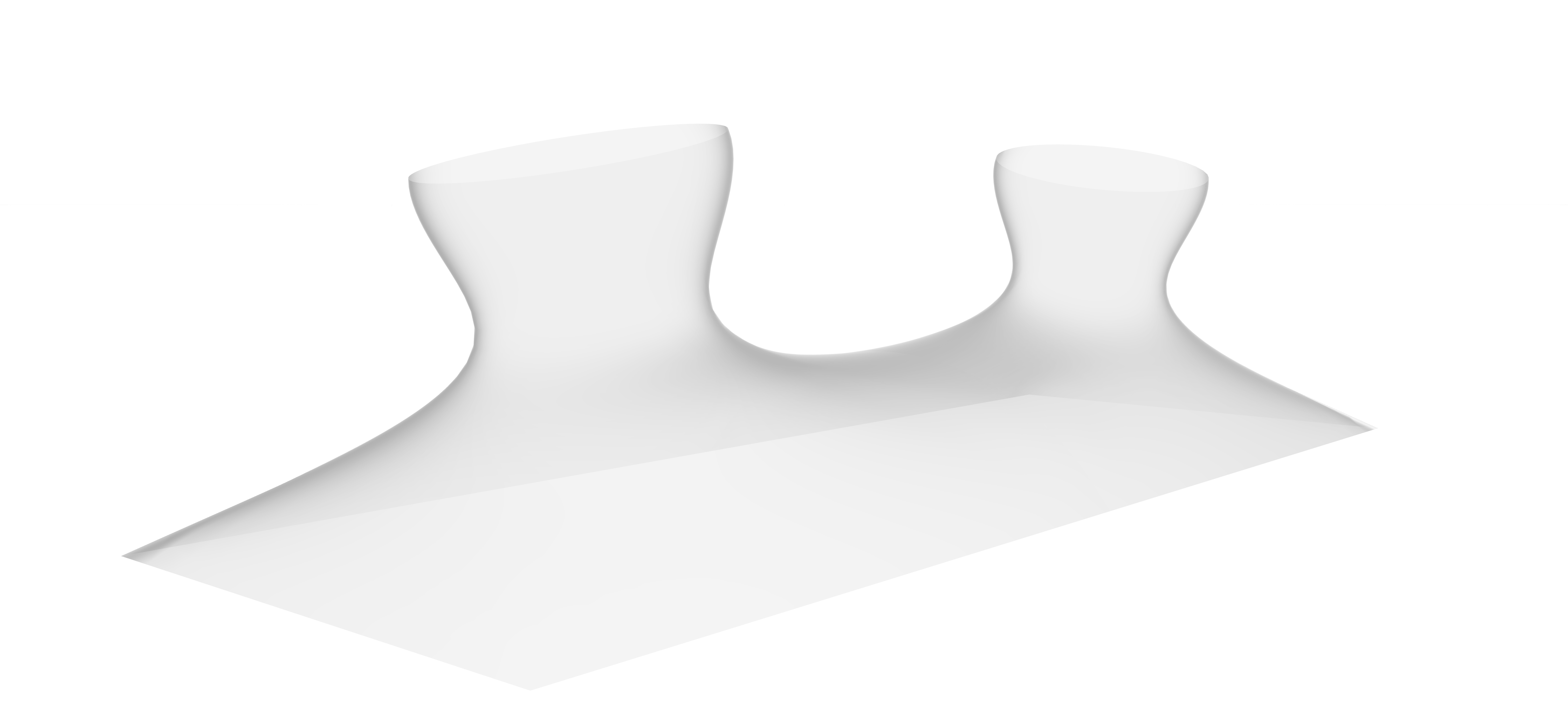}
    \caption{ It is possible to have boundary curves that lie in differently angled planes.}
    \label{fig:tiltedperspective}
\end{figure}

\begin{figure}[htp]
    \centering
    \begin{tikzpicture}
     \node (myfirstpic) at (0,0)   
    {\includegraphics[width=0.95\textwidth]{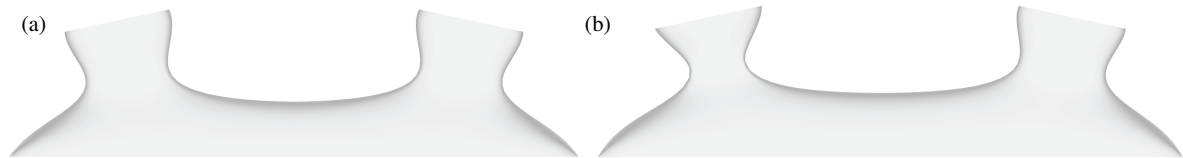}};
    \path (-7.5,1) node[below] {(a)};
    \path (0,1) node[below] {(b)};
    \end{tikzpicture}
    \caption{ Elevations of two surfaces sharing the same boundaries as in Fig. \ref{fig:tiltedperspective}. In (a) the currents in the two circular wires are equal, in (b) they are different.}
    \label{fig:titledcomparison}
\end{figure}

Figures \ref{fig:bridgeperspective} to \ref{fig:bridgsection} have a similar geometry to that in in Fig. \ref{fig:amancio}, but with several spans to model a potential bridge design. It is also possible to rotate the boundary curves as shown in \ref{fig:tiltedperspective} and \ref{fig:titledcomparison}.


\section{Conclusions and future work}
Constant solid angle surfaces enable one to control the boundary slope of a shell structure and hence achieve an approximately constant span-to-height ratio as the span varies. They also allow a principal curvature net to meet a plane boundary without cutting quadrilaterals. This means one can get a structurally viable shell also suitable for surface grids with planar panels.
\newline \indent
The constant boundary slope could also be used in the choice of Airy stress function where the constant slope would give a shell where the forces are concentrated at the corners of a boundary consisting of straight lines. This could be useful in projects with a similar context as the British Museum Great Court roof where forces are directed towards the corners relieving the walls from lateral thrust.   
\newline \indent
We have made no attempt to optimise the surfaces from the structural point of view, which depends both upon the shape of the surface and its boundary support \citep{GreenZerna68}. However, the conical shape at a boundary kink can be advantageous if there is a concentrated thrust at the kink, and that is the reason for the conical corners of the British Museum Great Court roof.
\newline \indent
The method could easily be adapted so that the required solid angle is no longer a constant, but some given function of the $x$, $y$, $z$ coordinates in space. Thus, if we were dissatisfied with the shape of the surfaces on a circular boundary in Fig. \ref{fig:circularboundary_1}, we could specify the solid angle as a function of $z$ or $\sqrt{x^2+y^2}$, which would preserve the rotational symmetry. In order to do this we would need to include the gradient of the required solid angle alongside the gradient in the solid angle.
\newline \indent
It is possible to generate surfaces having multiple complex boundary curves positioned and angled in different planes with individually tuned currents in the wires. One drawback is that the principal curvatures no longer follow the boundary curves. However, it is possible to rotate the slope around each curve independently by varying the current in the equivalent wire. Thus, the shapes can still be useful, and they can be used for other surface grids and possibly in other fields where one needs to generate surfaces without an initial mesh. However, there is still much to learn about the properties of constant solid angle surfaces. Having curves in different positioned and inclined planes as in Figs. \ref{fig:tiltedperspective} and \ref{fig:titledcomparison} it can be challenging to tune the parameters and find good initial positions for the points. Hence, further development of the technique can be done for such situations.

\section*{Acknowledgements}
We greatly appreciate the financial support from the Chalmers Foundation and the Digital Twin Cities Centre. We would also like to thank Daniel Piker for making us aware of the work regarding surfaces associated with knotted fields of by for instance \citet{Machon2013,Binysh2018, Binysh2019}. We are thankful for the helpful comments from Prof. Klas Modin. Furthermore,  we would like to thank Berlin Zoo for giving us the kind permission to use their picture from the Hippo House.

\bibliographystyle{IEEEtranN}
\bibliography{refs}
\end{document}